\documentclass[11pt]{article}

\renewcommand{\baselinestretch}{1.35}
\usepackage[marginratio=1:1,height=640pt,width=485pt,tmargin=72pt]{geometry}

\pdfoutput = 1

\usepackage{authblk}
\usepackage{amsmath}
\usepackage{amsthm}
\usepackage{bm}
\usepackage{amsfonts}
\usepackage{amssymb}
\usepackage{color}
\usepackage{graphicx}
\usepackage{natbib}
\usepackage{float}
\usepackage{wrapfig}
\usepackage{hyperref}
\usepackage{siunitx}
\usepackage{caption}
\usepackage{subcaption}
\usepackage{mathtools}
\usepackage[english]{babel}
\usepackage{mathrsfs}
\usepackage{longtable}
\allowdisplaybreaks
\usepackage{booktabs}
\usepackage{mathtools}
\setcitestyle{authoryear, round}
\usepackage[titletoc,title]{appendix}
\usepackage{algorithm2e} 
\usepackage{algorithmicx}
\RestyleAlgo{ruled}
\SetKwComment{Comment}{}{ }
\RestyleAlgo{ruled} 
\usepackage{rotating}
\usepackage{geometry}
\usepackage{booktabs}
\usepackage{threeparttable}
\usepackage[dvipsnames]{xcolor}
\usepackage{enumitem}
\usepackage{mathtools, nccmath, relsize}
\usepackage{amsthm}

\theoremstyle{plain}
\newtheorem{proposition}{Proposition}[section]

\theoremstyle{definition}

\theoremstyle{remark}
\newtheorem*{remark}{Remark}   

\allowdisplaybreaks
\providecommand{\keywords}[1]
{
  \small	
  \textbf{\textit{Keywords--}} #1
}

\title{Modular and Mobile Capacity Planning \\for Hyperconnected Supply Chain Networks}

\author[1,2]{Xiaoyue Liu}
\author[2,3]{Walid Klibi}
\author[1,2,4]{Benoit Montreuil}

\affil[1]{H.Milton Stewart School for Industrial \& Systems Engineering, Georgia Institute of Technology, Atlanta, United States}
\affil[2]{Physical Internet Center, Supply Chain \& Logistics Institute, Atlanta, United States}
\affil[3]{The Centre of Excellence for Supply Chain Innovation \& Transportation, Kedge Business School, Bordeaux, France}
\affil[4]{Coca-Cola Material Handling \& Distribution Chair, United States}

\date{\vspace{-0.35in} }

\begin{document} 
\maketitle

\begin{abstract}
The increased volatility of markets and the pressing need for resource sustainability are driving supply chains towards more agile, distributed, and dynamic designs. Motivated by the Physical Internet initiative, we introduce the Dynamic Stochastic Modular and Mobile Capacity Planning (DSMMCP) problem, which fosters hyperconnectivity through a network-of-networks architecture with modular and mobile capacities. The problem addresses both demand and supply uncertainties by incorporating short-term leasing of modular facilities and dynamic relocation of resources. We formulate DSMMCP as a partially adaptive multi-stage stochastic program that minimizes the expected multi-period costs under uncertainty. To tackle the inherent NP-hardness, we develop an enhanced stochastic dual dynamic integer programming (SDDiP) algorithm, which integrates strengthened cut generation, a tailored alternating cut strategy, and an efficient parallelization framework, and we establish structural dominance and monotonicity properties that formalize the value of the strengthened cuts and partial adaptivity. Numerical experiments inspired by a real case study of a large U.S. construction company demonstrate that the DSMMCP framework achieves approximately 15\% cost savings over static planning while improving resilience, reducing outsourcing costs, and supporting sustainability. Complementary experiments on synthetic instances confirm the effectiveness of the proposed SDDiP algorithm in terms of solution quality and runtime, as well as the scalability and robustness of the partially adaptive stochastic modeling framework across different network sizes and uncertainty levels.
\end{abstract}

\noindent \keywords{Modular capacity planning, Mobile production, Dynamic allocation and relocation, Multi-stage stochastic programming, Hyperconnected supply chains, Physical Internet}





\section{Introduction}
\label{sec:intro}

\subsection{Motivation}
\label{sec:mtv}

Volatility in resource availability places sustained operational pressure on capacity planning decisions, motivating the need for dynamic and adaptive capacity strategies \citep{vanmieghem2003capacity}. Growing resource scarcity has led to cost escalations, frequent supply disruptions, and decreased production capacities \citep{abb2024}, exposing the vulnerability of traditional supply chains and underscoring the need to invest sequentially in capacity to maintain resilience under disruptions \citep{klibi2024}. 

A persistent challenge is the under-utilization of fixed logistics assets, driven by mismatches between static capacity and volatile demand \citep{montreuil2011toward, Soundararajan2023, crainic2023hyperconnected}. Recent evidence illustrates these inefficiencies:
in Europe, average logistics vacancy exceeded 5\% in Q1 2025 \citep{cbre2025eu}, while U.S. industrial vacancy reached 6.6\% in Q2 2025, its highest level since 2014 \citep{cbre2025us}. Such idle assets create resource waste and hidden costs \citep{gill2008identifying} and limit the ability of firms to scale or adapt when conditions change, undermining both economic and resource sustainability. Volatile e-commerce volumes and short lead-time expectations further amplify uncertainty and increase the need for scalable, reconfigurable capacity \citep{niaz2022revolutionizing}. At the same time, digital technologies—such as real-time inventory visibility, predictive demand analytics, and platform-based fulfillment—raise expectations for speed and flexibility \citep{asawawibul2025influence, urban2025commerce} and drive firms toward more decentralized operational structures \citep{alarcon2022modular, kim2023network}.

These developments collectively motivate distributed network architectures whose capacity can be deployed closer to realized demand and adjusted dynamically over time. However, adopting such architectures is challenging, because establishing decentralized facilities typically requires substantial investments in infrastructure, equipment, and operations. To address this barrier, the concept of hyperconnected supply chain networks from the Physical Internet (PI) offers a promising alternative. The PI paradigm, introduced by \cite{montreuil2011toward}, promotes open asset sharing among certified users \citep{montreuil2013foundations}, enabling joint use of infrastructures, equipment, and services. Recent studies show that shared and mobile capacity can significantly improve utilization, responsiveness, and resilience \citep{kim2021hyperconnected, kulkarni2022resilient, faugere2022dynamic, liu2023logistics}.

\subsection{Problem Statement}
\label{sec:Prob stat}

This study is motivated by a collaborative project with a large modular construction company in the United States, which revealed the challenges of dynamically allocating modular and mobile capacities under volatile demand and uncertain supply availability.  
Similar operational questions arise across PI-related initiatives exploring shared and mobile infrastructure (e.g., \textit{MODULUSHCA}, \textit{CLUSTERS 2.0}, \textit{ICONET}).
Industry practice likewise reflects increased use of relocatable capacity—for example temporary fulfillment sites, dynamically repositioned parcel lockers, and logging camps in harvest planning \citep{LogiNext2025, faugere2022dynamic, jena2015modeling}.
In production contexts, \cite{marcotte2015modeling, marcotte2016introducing} introduced the concept of hyperconnected, modular, and mobile production and proposed optimization models for dynamically deploying production modules across facilities.  

\begin{figure}[h!]
    \centering
    \includegraphics[width=16cm]{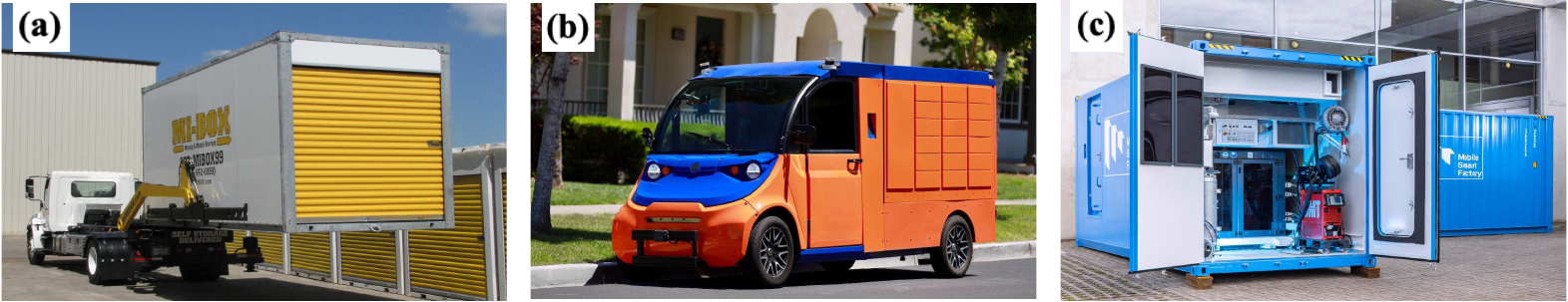}
    \caption{{Examples of modular and mobile capacities: (a) Mobile storage unit \citep{Mibox2024}, (b)
Mobile parcel lockers \citep{Symonds2019}, (c) Mobile production factory \citep{DefenceRedefined2023}.}}
    \label{fig:cap_examples}
    
\end{figure}

As illustrated in Figure \ref{fig:cap_examples}, modular and mobile capacities are no longer a distant vision but concrete solutions for hyperconnected supply chains, directly supporting the PI notion of an \textit{open network of networks}.  
From an operations research perspective, these developments motivate the formal definition of the Dynamic Stochastic Modular and Mobile Capacity Planning (DSMMCP) problem, which captures the tactical deployment of modular and mobile capacities in hyperconnected supply chains under both demand and supply uncertainty.

In contrast to classical facility-location models \citep{farahani2009facility}, DSMMCP incorporates short-term leasing, periodic redeployment, and stochastic demand and supply, yielding a dynamic NP-hard problem with temporal coupling.
While coalition-based approaches have highlighted the benefits of shared capacities \citep{roels2017win,jouida2021profit}, they face inherent scalability limitations as the number of participating actors increases. Hyperconnected networks instead support asset sharing at scale without explicit coalition formation, relying on standardized interfaces and open logistics infrastructures \citep{montreuil2013foundations, kim2021hyperconnected}.

To the best of our knowledge, no prior work provides a stochastic, multi-period optimization framework that captures the joint effects of mobility, modularity, and uncertainty in such networks. Addressing this gap is essential for quantifying the value of adaptive capacity strategies and for assessing when, where, and how such flexibility should be deployed in practice.

\subsection{Contribution}
\label{sec:ctb}

This paper introduces and studies a new tactical optimization problem—dynamic stochastic modular and mobile capacity planning (DSMMCP)—motivated by the emerging needs of hyperconnected supply chain networks. The main contributions are as follows.

\begin{itemize}

\item We frame and formalize the DSMMCP as a new multi-period stochastic problem for the tactical planning of modular and mobile capacities under demand and supply uncertainty.

\item We formulate DSMMCP as a partially adaptive multi-stage stochastic program (PAMSSP), and establish structural properties (e.g., monotone value of adaptivity) that justify the use of revision-point-based state relaxations and enable tractable solution approaches.

\item We develop an enhanced SDDiP algorithm tailored to PAMSSP, incorporating (i) strengthened Magnanti-Wong-type cuts, for which we establish dominance results, (ii) a novel multi-stage alternating cut strategy, and (iii) a parallelization framework for both forward and backward passes.

\item Using three years of real operational data from a large U.S. modular construction company, we demonstrate substantial improvements in cost, resilience, and resource utilization.

\end{itemize}

\vspace{0.01in}
The structure of the paper is as follows: Section \ref{sec:lr} reviews the literature; Section \ref{sec:problem} introduces DSMMCP and presents the PAMSSP formulation; Section \ref{sec:method} details the enhanced SDDiP algorithm; Section \ref{sec:ne} provides the case study and managerial insights; Section \ref{sec:con} concludes.

\section{Literature Review}
\label{sec:lr}

The classical Facility Location Problem (FLP), which centers on optimizing facility location decisions, has been extensively studied for decades \citep{armour1963heuristic, melo2009facility}. A major extension is the Dynamic Facility Location Problem (DFLP), which allows facility location changes over a multi-period planning horizon \citep{wesolowsky1973dynamic}. However, traditional DFLP models typically treat facilities as indivisible entities, with openings and closures occurring infrequently, which restricts network flexibility and leads to inefficient capital allocation.

\vspace{-0.01in}
{
\renewcommand{\baselinestretch}{0.9}\selectfont
\begin{table}[h!]
\centering
\caption{{Problem decisions and settings in studies related to the Dynamic Facility Location Problem with Modular Capacities.}}
\label{tab:lr_problem}
\resizebox{0.93\textwidth}{!}{
\begin{threeparttable}
  \begin{tabular}{lccccccl}\hline
    & \multicolumn{3}{c}{\textbf{Uncertainty}} & \multicolumn{2}{c}{\textbf{Decisions}} & \multicolumn{2}{c}{\textbf{Problem settings}} \\\cmidrule(lr){2-4} \cmidrule(lr){5-6} \cmidrule(lr){7-8} 
    \textbf{Reference} & \textbf{D/S} \tnote{a} & \textbf{Demand} & \textbf{Supply} & \textbf{O/C/R} \tnote{b} & \textbf{E/R/RL} \tnote{c} & \textbf{TD} \tnote{d} & \textbf{Context} \tnote{d} \\ \hline
    \cite{jena2015dynamic, jena2017lagrangian} & D &  &  & O/C/R & E/R & & N/S \\
    \cite{becker2019value} & D &  &  & O/-/- & E/R/RL &  & Chemical industry \\ 
    \cite{hong2020optimal} & D &  &  & O/C/R & E/R/RL &  & Gas industry \\ 
    \cite{allman2020dynamic} & D &  &  & O/-/- & E/RL &  & Chemical industry \\ 
    \cite{correia2021integrated} & S & \checkmark &  & O/-/- & E/R & \checkmark & N/S \\ 
    \cite{faugere2022dynamic} & S & \checkmark &  & -/-/- & RL & \checkmark & Urban parcel delivery \\ 
    \cite{kayser2023relocatable} & S & \checkmark &  & O/C/R & E/R/RL &  & N/S \\ 
    \cite{ge2023multistage} & S & \checkmark &  & -/-/- & E/RL &  & Closed-loop supply chain \\ 
    \cite{correia2024matheuristic} & D &  &  & O/C/R & E/R &  & N/S \\
    \cite{vstadlerova2024solving}& S & \checkmark &  & O/C/R & E/R &  & Hydrogen industry \\ 
    \cite{liu2025dynamic} & S & \checkmark &  & O/C/R & E/R & \checkmark & Freight transportation \\ 
    Our work & S & \checkmark & \checkmark & O/C/R & E/R/RL & \checkmark & N/S \\\hline
  \end{tabular}
  \begin{tablenotes}
    \footnotesize
    \item[a] D: Deterministic; S: Stochastic.
    \item[b] O/C/R: The status of facilities can be opening (O), closing (C), or reopening (R). Note: "-" indicates that there is no associated cost term for changing the facility from a previous state to that specific state.
    \item[c] E/R/RL: The number of capacity modules can be expanded (E) and reduced (R), and we allow capacity modules to be relocated (RL). 
    \item[d] TD: Time decoupling, where the problem has multiple time scales; N/S: Not specified (i.e., generalized context). 
  \end{tablenotes}
\end{threeparttable}
}
\end{table}
}

To address these limitations, the Dynamic Facility Location Problem with Modular Capacities (DFLPM) has emerged as a more flexible framework. As summarized in Table~\ref{tab:lr_problem}, recent studies consider incremental capacity adjustments through facility opening, closing, and reopening  (referred to as O/C/R), as well as capacity expansion and reduction (referred to as E/R) via modular units \citep{jena2015dynamic, correia2024matheuristic}. This modular structure shifts planning from a purely strategic level toward tactical or operational decision-making \citep{ge2023multistage}, while substantially increasing modeling and computational complexity.

A more advanced research stream further incorporates mobility by allowing physical relocation of capacity modules within the network (referred to as RL in Table~\ref{tab:lr_problem}). Such mobility improves resource utilization and sustainability by redeploying existing capacity to meet spatially and temporally shifting demand. Applications include relocatable processing devices in the energy sector \citep{hong2020optimal} and mobile logistics assets in urban delivery systems \citep{faugere2022dynamic}. However, very few studies have combined modular capacity planning with full flexibility and mobility. Even fewer have applied this to hyperconnected supply chain networks, where the sharing of logistics capacity and infrastructure among network participants poses a major planning challenge while offering substantial potential for efficiency and sustainability.

Another important modeling feature is time decoupling (referred to as TD in Table~\ref{tab:lr_problem}), where strategic or tactical network decisions are made on a coarser time scale than operational decisions \citep{correia2021integrated, liu2025dynamic}. This study adopts this feature through a partially adaptive approach, as detailed in Sections~\ref{sec:pd} and~\ref{sec:PAMSSP}, to better reflect real-world planning processes.

Uncertainty is inherent in dynamic facility location and capacity planning problems, motivating stochastic modeling approaches. As shown in Table~\ref{tab:lr_problem}, existing stochastic DFLPM studies primarily focus on demand uncertainty. Several works incorporate stochastic demand to guide facility openings and capacity adjustments, mitigating the cost of frequent network reconfiguration \citep{kayser2023relocatable, vstadlerova2024solving}. In contrast, supply-side uncertainty has received limited attention. A key contribution of this study is the integration of both demand and supply uncertainty, accounting for disruptions that affect the effective availability of modular capacity.

{
\renewcommand{\baselinestretch}{0.9}\selectfont
\begin{table}[h!]
\centering
\caption{{Modeling and solution approach in studies related to Stochastic Dynamic Facility Location Problems.}}
\label{tab:lr_modeling_solution}
\resizebox{\textwidth}{!}{
\begin{threeparttable}
  \begin{tabular}{lcccccccc}\hline
    & \multicolumn{3}{c}{\textbf{Modeling}} & \multicolumn{2}{c}{\textbf{Solution approach}} & \multicolumn{3}{c}{\textbf{Numerical experiments}} \\\cmidrule(lr){2-4} \cmidrule(lr){5-6} \cmidrule(lr){7-9} 
    \textbf{Reference} & \textbf{Model} \tnote{a} & \textbf{I/FCC} \tnote{b} & \textbf{DS/O} \tnote{b} & \textbf{Method} \tnote{c} & \textbf{Detail }\tnote{d} & \textbf{Case study}  & \textbf{Horizon} \tnote{c} & \textbf{\# Periods} \\ \hline
    \cite{correia2018stochastic} & TSSP & \checkmark &  & E & CPLEX; VI & & N/S & 3 \\ 
    \cite{correia2021integrated} & TSSP & \checkmark &  & E & CPLEX; VI & & N/S & 24 \\ 
    \cite{arslan2021distribution} & TSSP & \checkmark & \checkmark & E & L-SH & \checkmark & 1 month & 30 \\ 
    \cite{faugere2022dynamic} & TSSP & \checkmark & \checkmark & E & BD; RHA & \checkmark & 2 months & 8 \\  
    \cite{kayser2023relocatable} & TSSP & \checkmark & \checkmark & E & CPLEX & & N/S & 12 \\ 
    \cite{ge2023multistage} & MSSP & \checkmark &  & E & SDDiP & & N/S & 8 \\ 
    \cite{yang2023distributionally} & DRO &  & \checkmark & E & B\&BC & \checkmark & 5 days & 5 \\ 
    \cite{vstadlerova2024solving} & MSSP & \checkmark & \checkmark & H & LH & \checkmark & 2 weeks & 14 \\
    \cite{wang2024global} & TSSP & \checkmark &  & E & BD &  & 1 year & 24 \\ 
    \cite{liu2025dynamic} & TSSP & \checkmark & \checkmark & E & B\&BC & \checkmark & 1 year & 52 \\ 
    Our work & PAMSSP & \checkmark & \checkmark & E & SDDiP & \checkmark & 1 year & 12 \\ \hline
  \end{tabular}
  \begin{tablenotes}
    \footnotesize
    \item[a] TSSP: Two-stage stochastic program; MDP: Markov decision process; IP: Integer program; MOP: Multi-objective program; MSSP: Multi-stage stochastic program; DRO: Distributionally robust optimization; PAMSSP: Partially adaptive multi-stage stochastic program. 
    \item[b] I/FCC: Inventory/flow conservation constraints; DS/O: Demand selection/outsourcing.
    \item[c] H: Heuristics; E: Exact method; N/S: Not specified.
    \item[d] LH: Lagrangian-based heuristic; VI: valid inequalities; L-SH: L-shaped method; BD: Benders decomposition; RHA: Rolling-horizon approach; DDiP: Stochastic dual dynamic integer programming; B\&BC: Branch-and-Benders-cut algorithm.
  \end{tablenotes}
\end{threeparttable}
}
\end{table}
}

The dual challenges of dynamics and stochasticity inherent in this problem necessitate efficient modeling and solution approaches. We thus turn our review to methods developed for the broader class of Stochastic Dynamic Facility Location Problems (SDFLP), as summarized in Table \ref{tab:lr_modeling_solution}. Most studies adopt two-stage stochastic programs (TSSP), which offer computational tractability but restrict adaptivity by fixing strategic decisions before uncertainty is revealed \citep{correia2018stochastic}. Multi-stage stochastic programs (MSSP) allow decisions to evolve as information is revealed, but are computationally demanding \citep{ge2023multistage}. To balance this trade-off, \cite{zou2018partially} proposed a partially adaptive multi-stage stochastic program (PAMSSP), which allows state decisions to be adaptive up to certain prescribed nodes before restricting them to a two-stage structure. Building on this, \cite{basciftci2024adaptive, kayacik2025partially} developed a generalized framework by allowing preset (or optimize) when state variables can be revised within the planning horizon. As far as we know, our work presents the first application of PAMSSP to the dynamic stochastic modular and mobile capacity planning problem, demonstrated through a real case study.

Moreover, we note two modeling features critical to our problem. The first is the use of inventory or flow conservation constraints (referred to as I/FCC in Table \ref{tab:lr_modeling_solution}), which link decisions across time periods and are essential for capturing system dynamics \citep{wang2024global}. The second is the inclusion of demand selection or outsourcing (referred to as DS/O in Table \ref{tab:lr_modeling_solution}) to provide crucial operational flexibility in managing capacity shortfalls \citep{arslan2021distribution, yang2023distributionally}.

Since the FLP is NP-hard, its dynamic and stochastic variants pose significant computational challenges. As shown in Table~\ref{tab:lr_modeling_solution}, existing solution approaches include both heuristics and exact methods, with exact methods predominantly relying on decomposition techniques designed to break the problem down into more manageable subproblems. For multi-stage stochastic mixed-integer programs, SDDiP has emerged as a powerful solution framework. Early work introduced nested Benders decomposition for multi-stage stochastic programming and led to the stochastic dual dynamic programming (SDDP) method \citep{birge1985decomposition, pereira1991multi}. While these approaches focused on convex problems, the SDDiP algorithm proposed by \cite{zou2019stochastic, zou2018multistage} extends this framework to multi-stage stochastic integer programs. Subsequent studies demonstrated its ability to handle mixed-integer state variables and large-scale scenario trees \citep{lara2020electric, quezada2022combining}. Given the structural similarities between PAMSSP and MSSP, we develop a tailored SDDiP-based algorithm to solve the proposed PAMSSP formulation, representing the first application of SDDiP in this modeling class.

\section{Dynamic Stochastic Modular and Mobile Capacity Planning}
\label{sec:problem}

\subsection{Problem Description}
\label{sec:pd}

Building on the problem statement in Section~\ref{sec:Prob stat}, we now provide a formal description of the \textit{Dynamic Stochastic Modular and Mobile Capacity Planning (DSMMCP)} problem. The setting is motivated by a company that operates a hyperconnected supply chain network, relying on short-term leases and the dynamic deployment of modular and mobile capacities to meet customer demand. Two main sources of uncertainty are considered: \textit{demand} and \textit{supply}. At the beginning of the planning horizon, the set of potential customer locations is known, but demand volumes are uncertain and revealed progressively over time. On the supply side, while each capacity module has a fixed nominal throughput, its actual availability may be reduced by disruptions such as equipment failures, labor shortages, or hazard disruptions.  

\begin{figure}[h!]
    \centering
    \includegraphics[width=15cm]{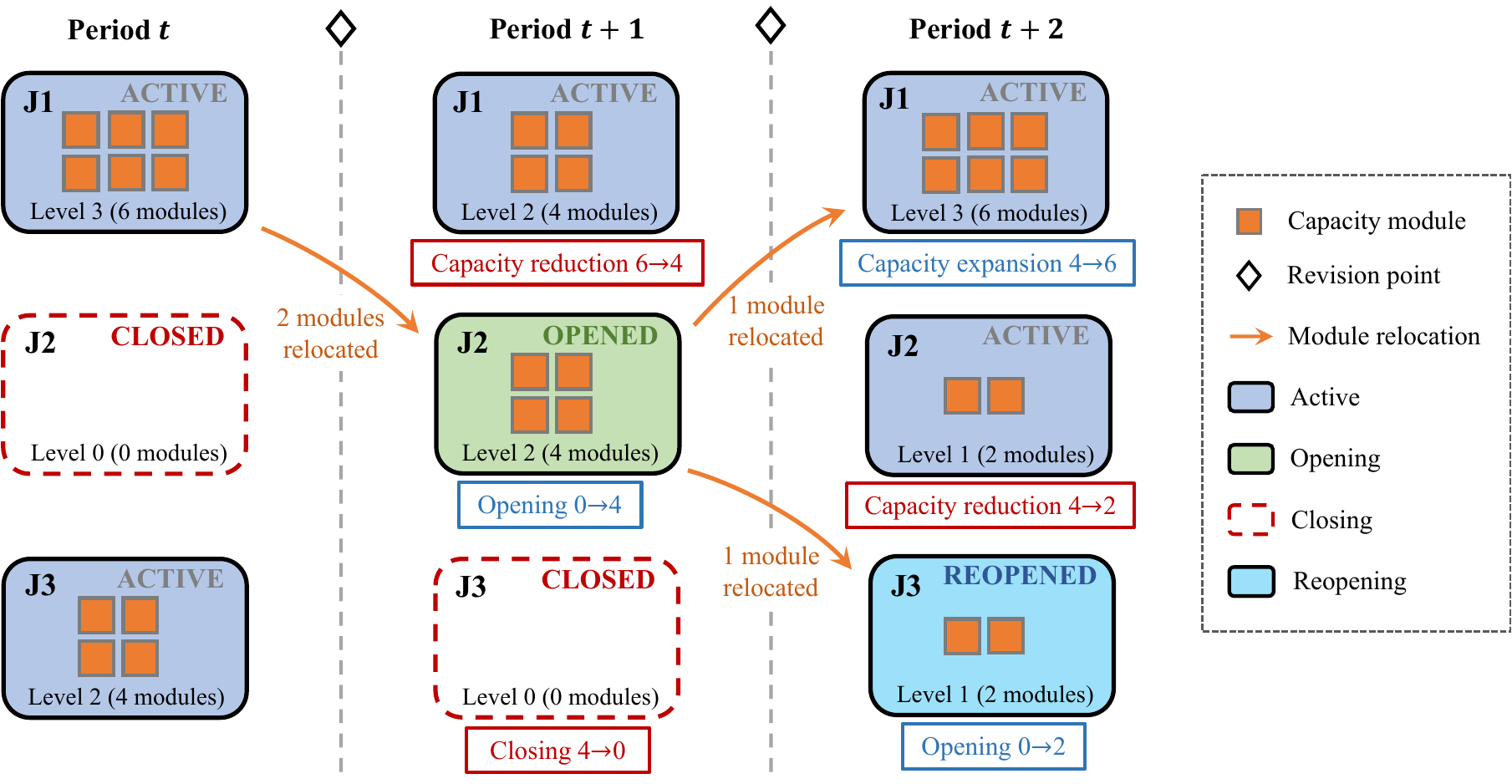}
    \caption{\small{Dynamic modular and mobile capacity planning across three planning periods with two revision points.}}
    \label{fig:relocation}
    \vspace{-0.1in}
\end{figure} 

To mitigate these uncertainties, the company can revise its network quarterly as new information becomes available. Temporary facilities may be leased near demand locations to reduce shipping distances and service delays, while capacity modules can be reallocated among facilities to better align with demand. By dynamically adjusting modular and mobile capacities, facilities can flexibly scale up or down, open, close, or reopen as needed, allowing the network to adapt to evolving demand and absorb supply disruptions. A facility can expand capacity by renting new modules from the depot or by relocating modules from other facilities. Conversely, capacity can be reduced by returning unused modules to the depot or reallocating them elsewhere. Depending on the context, the depot may represent a lessor, a merchant, or another type of capacity provider. Figure \ref{fig:relocation} demonstrates key operational capabilities in our capacity planning framework: facility opening, closing, and reopening; capacity expansion and reduction; module relocation between facilities (orange arrows). In line with practice, we assume that potential facility locations are predetermined, reflecting prior screening for infrastructure availability or strategic suitability.  
In addition, the company may outsource excess demand when assigned capacity is insufficient. Outsourcing, however, incurs a penalty cost that reflects either payments to external providers or the use of costly on-demand capacity. A higher penalty incentivizes greater in-house deployment, whereas a lower penalty offers flexibility by allowing leaner network configurations and increased outsourcing.   

\begin{figure}[h!]
    \vspace{-0.01in}
    \centering
    \includegraphics[width=16.5cm]{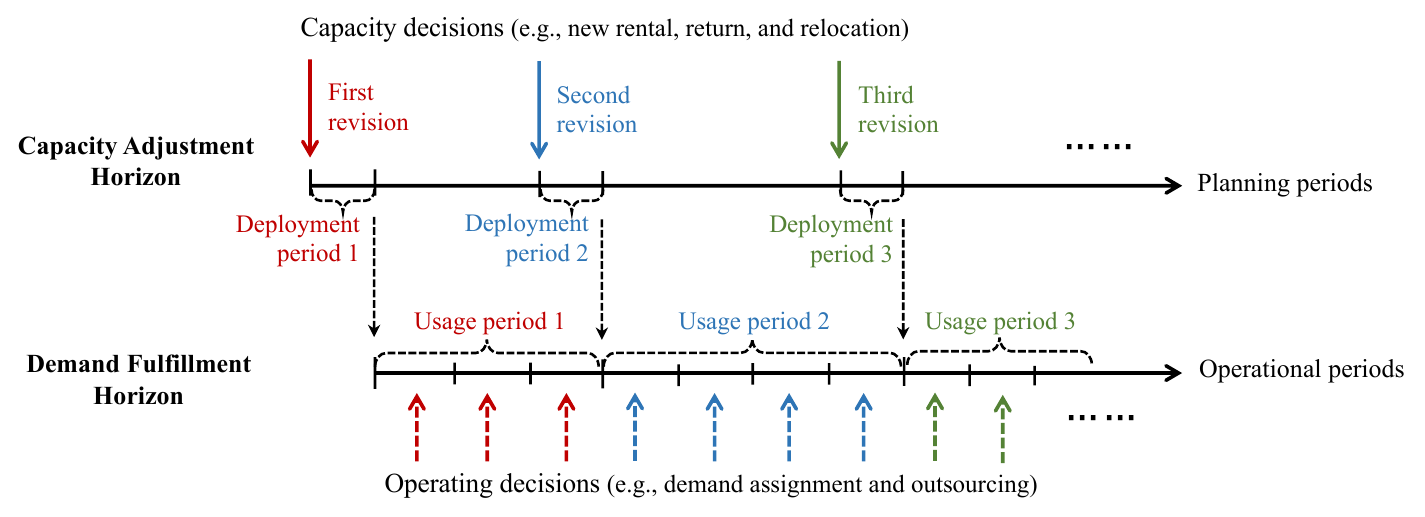}
    \caption{{Timeline of the multi-stage decision-making process in DSMMCP.}}
    \label{fig:timeline}
    \vspace{-0.05in}
\end{figure} 

As illustrated in Figure~\ref{fig:timeline}, the resulting problem is framed as a multi-stage decision process and features the time decoupling (TD) with two distinct horizons: a capacity adjustment horizon, where tactical capacity decisions are made periodically to reflect common contractual constraints (e.g., minimum lease terms), and a demand fulfillment horizon, where operating decisions occur at a much finer time scale (i.e., multi-horizons, \citep{kaut2014multi}). Specifically, at the beginning of the horizon, the company makes capacity deployment decisions (red solid arrow). As time progresses and uncertainty is revealed, these decisions may be revised at the beginning of each subsequent period (green and blue arrows). In addition, each tactical capacity decision is followed by a deployment lag, representing the implementation time before the capacity revision becomes available. After the deployment period, operational decisions such as demand assignment and outsourcing are made during the corresponding usage period, given the current capacity and realized uncertainties (dashed arrows). 

This multi-stage structure enables more effective capital management by replacing a substantial upfront expenditure with a series of phased investments that can be dynamically adjusted in response to market growth and evolving system parameters \citep{singh2009dantzig, huang2009value, jaoui2025multi}. Moreover, the deployment lag depicts the information asymmetry between different decision types. Tactical decisions are determined proactively based on the information available at the revision point, whereas operational decisions are made reactively, benefiting from more information revealed by the beginning of the usage period. While providing realism and flexibility, this multi-stage structure, combined with underlying uncertainties, suffers from the well-known "curse of dimensionality" \citep{pereira1991multi}. The number of states and future scenarios grows exponentially with the number of potential facilities, capacity levels, and planning periods, making the problem computationally intractable to solve optimally using exact methods.

\subsection{Scenario-based Uncertainty Modeling}
\label{sec:str}
Due to the presence of both demand and supply uncertainties, the related information is unavailable during the tactical planning phase but plays a key role in shaping the capacity decisions. To capture this bottom-up influence and properly model these uncertainties, we adopt a scenario tree approach, modeling all uncertain parameters as discrete-time stochastic processes with a finite set of realizations. Let $\boldsymbol{T} = \{1,\dots,t,\dots,T\}$ represent the planning horizon, where each implementation period $t$ corresponds to a month. The decision-making process unfolds over $T$ stages, with each period denoting a single stage. 

\begin{figure*}[h!]
\centering
\begin{minipage}[t]{0.4\textwidth}

\centering
\includegraphics[width=0.9\linewidth]{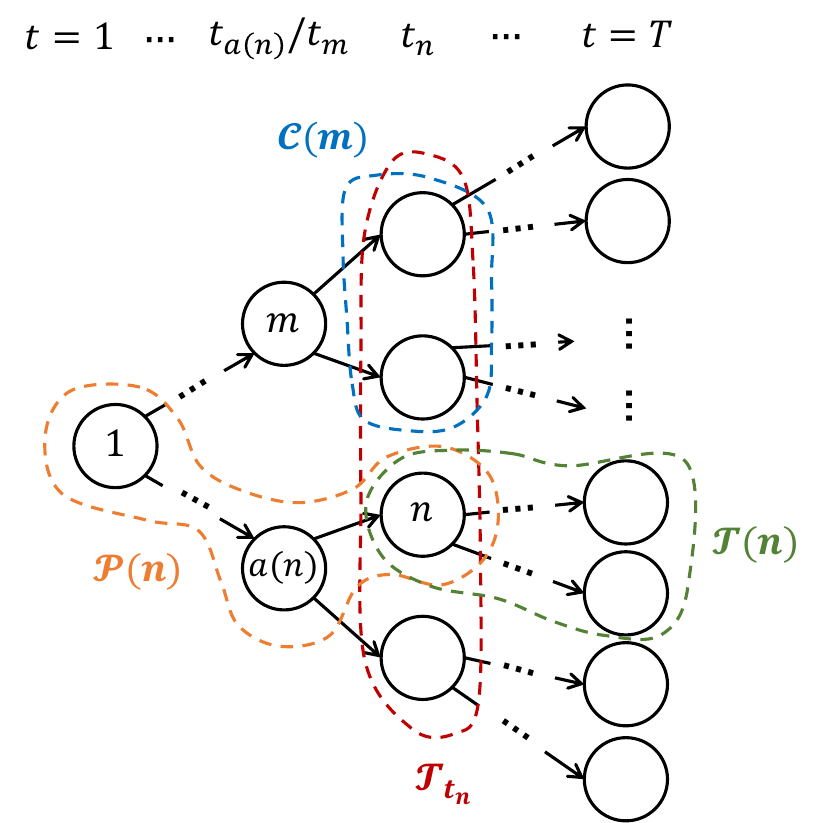}
\captionof{figure}{Illustration of scenario tree $\boldsymbol{\mathcal{T}}$.}
\label{fig:scenario}
\end{minipage}\hfill
\begin{minipage}[t]{0.52\textwidth}

\centering
\includegraphics[width=0.9\linewidth]{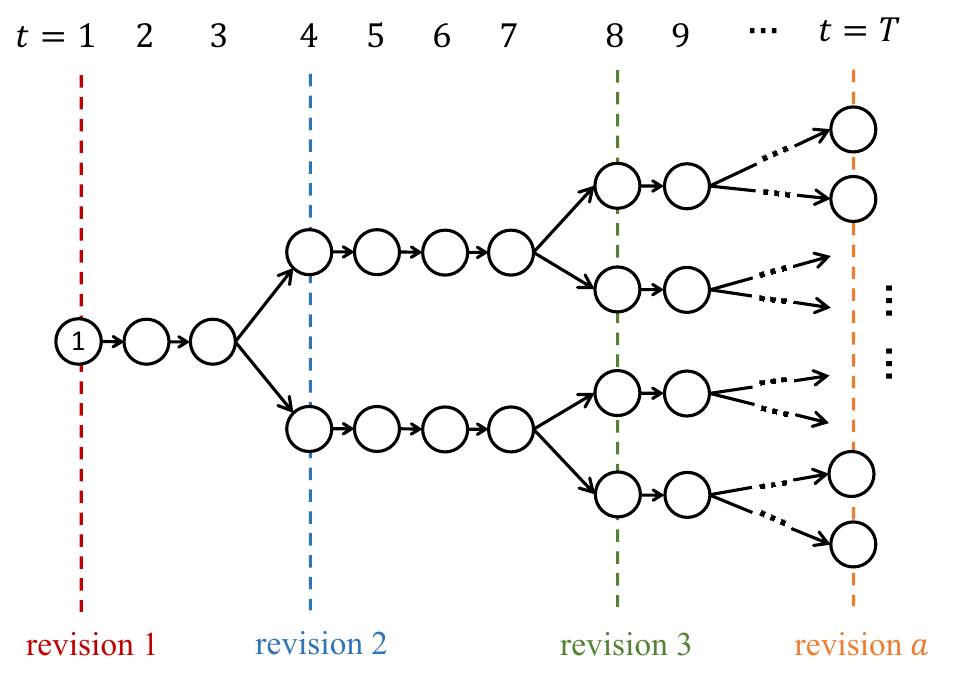}
\captionof{figure}{Partially adaptive decision structure of state variables with $a$ revisions.}
\label{fig:scenario_adaptive}
\end{minipage}
\end{figure*}

Figure \ref{fig:scenario} shows the structure of the scenario tree $\boldsymbol{\mathcal{T}}$, with each node representing a realization of the underlying uncertainties. Except for the root node (i.e., $n=1$), each node $n$ in stage $t$ has a unique ancestor node in the previous stage $t-1$, denoted by $a(n)$. One or multiple child nodes in stage $t+1$ are connected with node $n$, and the set of these child nodes is represented as $\boldsymbol{\mathcal{C}(n)}$. We also define the period of node $n$ as $t_n$ and the stage including node $n$ as $t_n$. The set of all nodes in the same stage as node $n$ is denoted by $\boldsymbol{T(n)}$, while $\boldsymbol{\mathcal{T}_t}$ represents the set of all nodes in stage $t$. In addition, let $\boldsymbol{\Omega_t}$ represent the set of uncertain parameter realizations up to period $t$. We denote $\boldsymbol{\omega_{tm}}$ as the $m$-th realization of $\boldsymbol{\Omega_t}$ up to period $t$ and $N_{t}$ as the total number of realizations in stage $t$. The complete sequence of realizations from period $1$ to $T$ forms a scenario. To simplify the notation, we define the set of scenarios as $\boldsymbol{\Omega}$, meaning we have $\boldsymbol{\Omega_T} \equiv \boldsymbol{\Omega}$ and $\boldsymbol{\omega_{Tm}} \equiv \boldsymbol{\omega_{m}}, \forall m$. Furthermore, the unique realization from the root node $0$ to node $n$ is expressed as $\boldsymbol{\omega(n)}$. If node $n$ is a leaf node (i.e., $n \in \boldsymbol{\mathcal{T}_T}$), realization $\boldsymbol{\omega(n)}$ denotes scenario $\boldsymbol{\omega}$. If node $n$ is a non-leaf node, $\boldsymbol{\mathcal{T}(n)}$ is the subtree rooted at node $n$. Lastly, the probability of node $n$ is given by $p_n$, which is the probability of the realization sequence $\boldsymbol{\omega_{t_n}}$. The summation of probabilities of all nodes within one stage is 1, and the summation of probabilities of all child nodes is equal to the probability of the ancestor node. 

\vspace{-0.1in}
\subsection{Partially Adaptive Multi-Stage Stochastic Program (PAMSSP)}
\label{sec:PAMSSP}

There are mainly two types of decision variables: state variables (e.g., capacity decisions), which link different stages, and local variables (e.g., operating decisions), which are specific to node $n$'s problem. In an MSSP, both variable types can be adjusted at each stage based on observed uncertainties, enabling full adaptivity but at an extremely high computational cost. As an alternative, the TSSP approximation is often utilized to simplify MSSP by enforcing equality in state variables at each stage $t$. This restriction, however, sacrifices the flexibility to adapt to uncertainty over time. To achieve a balance between decision adaptability and computational complexity, in this study, we formulate DSMMCP as a PAMSSP. We remark that the revision timing decisions in PAMSSP can be treated in two ways: either as input parameters predefined by the decision maker or as here-and-now decisions made before uncertainties unfold. While state variables can only be revised at these revision points, local variables can be updated at every stage. In our problem, this implies that the timing of capacity adjustments is determined at the beginning of the planning horizon, while the magnitude of adjustments is decided based on the information available at each decision point. This setting reflects real-world practice, where companies are often required to predefine the duration of leasing contracts with facility providers. To maintain contract stability, we limit flexibility by restricting at most $a$ revision points, where $a \in \{1, ..., T\}$. Setting $a=1$ reduces PAMSSP to a traditional TSSP, while $a=T$ makes PAMSSP equivalent to an MSSP. 

Figure \ref{fig:scenario_adaptive} illustrates a decision structure of state variables in a PAMSSP. In the figure, stages 4 and 8 are two predetermined revision points, allowing the state variables at nodes in stage 4 (respectively, stage 8) to adjust based on uncertainties observed up to stage 4 (respectively, stage 8). However, within each sub-tree rooted at a stage 4 node, the nodes from stages 4 to 7 are compressed, indicating a single state decision is applied across these nodes. In our context, this represents that the company updates the facility capacity decisions at the beginning of the fourth month and commits to these decisions for the following three months, thereby signing a 4-month leasing contract with the facility provider.
\subsubsection*{Value of Partial Adaptivity}

The partially adaptive multistage stochastic program (PAMSSP) interpolates between the
two-stage static model (TSSP), where no recourse in the state variable is allowed, and the
fully adaptive multistage model (MSSP), where the state may branch freely across the
scenario tree. Increasing the number of revision points expands the feasible set by
relaxing nonanticipativity constraints on the capacity-state variable~$Y$. Therefore, the
optimal expected cost cannot increase when additional revision points are introduced.


\begin{proposition}[Monotone Value of Partial Adaptivity]\label{prop:adaptivity}
Let $z^{\mathsf{MSSP}}$ denote the optimal value of the fully adaptive multistage model,
$z^{\mathsf{TSSP}}$ the optimal value of the two-stage static model,
and $z^{\mathsf{PAMSSP}}(a)$ the optimal value of the partially adaptive model with at
most $a$ revision points. Then for any integers $1 \le a_1 \le a_2 \le T$,

{
\[
z^{\mathsf{MSSP}}
\;\le\;
z^{\mathsf{PAMSSP}}(a_2)
\;\le\;
z^{\mathsf{PAMSSP}}(a_1)
\;\le\;
z^{\mathsf{TSSP}},
\]
and equivalently,
\[
z^{\mathsf{PAMSSP}}(a) - z^{\mathsf{PAMSSP}}(a+1) \ge 0 \quad \forall\, a < T.
\]
}
\end{proposition}

This value-of-information property shows that additional revision points can only improve
(or at worst leave unchanged) the expected cost. In our computational results, we observe
that a small number of well-placed revision points already captures most of the value of
full adaptivity. A complete proof is provided in
Appendix~\ref{appendix:pamssp_proof}.

\subsection{Mathematical Formulation}
\label{sec:formulation} 

Using a scenario tree to represent the evolution of the uncertain parameters, we formulate the problem as a PAMSSP, which can be solved using commercial optimization solvers. The notation is listed in Appendix \ref{appendix_notation}.

\begin{flalign}
    \min \ & 
    \sum_{n \in \boldsymbol{\mathcal{T}}} p_n  \left(
    \sum_{j \in \boldsymbol{J}} \sum_{l_1 \in \boldsymbol{L}} \sum_{l_2 \in \boldsymbol{L}} c_{j l_1 l_2 n}^F  Y_{j l_1 l_2 n}
    +\sum_{(j,j^{'}) \in \boldsymbol{\mathit{JJ}}} c_{j j^{'} n}^A  F_{jj^{'}n} 
    + \sum_{(i,j) \in \boldsymbol{\mathit{IJ}}} c_{i j n}^T  X_{ijn} 
    + \sum_{j \in \boldsymbol{J}} c_{n}^O  R_{jn}\right) \hspace{-0.1in} & \label{eq:1}
\end{flalign}
s.t.,
\begin{flalign}
    & \sum_{j \in \boldsymbol{J} \mid (i,j) \in \boldsymbol{\mathit{IJ}}} X_{ijn} \geq d_{in} 
    & \forall i \in \boldsymbol{I}, n \in \boldsymbol{\mathcal{T}} \label{eq:2} \\
    & k_{jn} S_{jn} + R_{jn} \geq \sum_{i \in \boldsymbol{I} \mid (i,j) \in \boldsymbol{\mathit{IJ}}} X_{ijn} 
    & \forall j \in \boldsymbol{J}, n \in \boldsymbol{\mathcal{T}} \label{eq:4} \\
    & S_{j1} = u_{j v_j} + \sum_{j^{'} \in \boldsymbol{J} \cup \{0\} \mid (j,j^{'}) \in \boldsymbol{\mathit{JJ}}} \left(F_{j^{'}j1} - F_{jj^{'}1}\right) 
    & \forall j \in \boldsymbol{J} \label{eq:4.5} \\
    & S_{jn} = S_{j, a(n)} + \sum_{j^{'} \in \boldsymbol{J} \cup \{0\} \mid (j,j^{'}) \in \boldsymbol{\mathit{JJ}}} \left(F_{j^{'}jn} - F_{jj^{'}n}\right) 
    & \forall j \in \boldsymbol{J}, n \in \boldsymbol{\mathcal{T}} \backslash\{1\} \label{eq:5} \\
    & S_{jn} = \sum_{l_1 \in \boldsymbol{L}} \sum_{l_2 \in \boldsymbol{L}} u_{j l_2}  Y_{j l_1 l_2 n} 
    & \forall j \in \boldsymbol{J}, n \in \boldsymbol{\mathcal{T}} \label{eq:6}  \\
    & \sum_{l \in \boldsymbol{L}} Y_{j v_{j} l 1}=1 
    & \forall j \in \boldsymbol{J} \label{eq:7}\\
    & \sum_{l_1 \in \boldsymbol{L}} Y_{j, l_1, l, a(n)}=\sum_{l_2 \in \boldsymbol{L}} Y_{j, l, l_2, n} 
    & \forall j \in \boldsymbol{J}, l \in \boldsymbol{L}, n \in \boldsymbol{\mathcal{T}} \backslash\{1\} \label{eq:8} \\
    & \sum_{l_1 \in \boldsymbol{L}} \sum_{l_2 \in \boldsymbol{L}} Y_{j l_1 l_2 n}=1 
    & \forall j \in \boldsymbol{J}, n \in \boldsymbol{\mathcal{T}} \label{eq:9}
\end{flalign}

\begin{flalign}
    & Y_{j l_1 l_2 n} \leq w_{t}
    & \forall j \in \boldsymbol{J}, l1, l2 \in \boldsymbol{L} \text{ if } l1 \neq l2, n \in \boldsymbol{\mathcal{T}_t}, t \in \boldsymbol{T} \label{eq:10}\\
    & Y_{j l_1 l_2 n}  \in \{0, 1\} 
    & j \in \boldsymbol{J}, l_1, l_2 \in \boldsymbol{L}, n \in \boldsymbol{\mathcal{T}} \label{eq:15} \\
    & F_{jj^{'}n}, S_{jn} \in \mathbb Z_{\ge 0} 
    & \forall (j,j^{'}) \in \boldsymbol{\mathit{JJ}}, j \in \boldsymbol{J}, n \in \boldsymbol{\mathcal{T}} \label{eq:16} \\
    & X_{ijn}, R_{jn} \geq 0 
    & \forall (i,j) \in \boldsymbol{\mathit{IJ}}, j \in \boldsymbol{J}, n \in \boldsymbol{\mathcal{T}} \label{eq:17}.
\end{flalign}

In the objective function (\ref{eq:1}), we aim to minimize total costs, including facility capacity adjustment costs, capacity module deployment costs, transportation costs, and outsourcing costs. Constraints (\ref{eq:2}) state that the company is required to serve all customer demand in every period. Constraints (\ref{eq:4}) indicate that each facility has two options to finish its assigned demand: either use its own production capacity or outsource part of the demand using external capacity acquisition. 
Flow conservation constraints (\ref{eq:4.5})-(\ref{eq:5}) enforce that the movement of capacity modules is balanced with the number available at each facility over time, while constraints (\ref{eq:6}) define the relationship between a facility’s capacity and the number of deployed capacity modules. Constraints (\ref{eq:7})-(\ref{eq:8}) specify the initial capacity level of each facility and ensure that changes in capacity levels are consecutive from one period to another. Constraints (\ref{eq:9}) guarantee that exactly one capacity level is selected for each facility in each period. Constraints (\ref{eq:10}) ensure that capacity decisions can only be revised at predetermined revision points, representing the decision structure for partially adaptive modeling. More specifically, $w_t$ is a binary parameter that indicates whether the capacity revision is allowed at period $t$. If $w_t=1$ (i.e., we are at the revision point), constraints (\ref{eq:10}) become $Y_{j l_1 l_2 n} \leq 1$, so there is no restriction on capacity revision at the node $n$. If $w_t=0$ (i.e., we are not at the revision point), constraints (\ref{eq:10}) force $Y_{jl_1l_2n}=0$ for these $l1 \neq l2$. Then, together with constraints (\ref{eq:8})-(\ref{eq:9}), we have $Y_{jl l n}=1 \forall j, n$, where $l$ is the same capacity level as the parent node $a(n)$. This means that all nodes in this period keep the same capacity level as their last period (i.e., non-anticipativity is enforced).

\section{SDDiP Framework}
\label{sec:method}

\subsection{Reformulation Based on Dynamic Programming and Decomposition Rationale}
\label{sec:dp}
The proposed formulation (\ref{eq:1})-(\ref{eq:17}) exhibits an approximate block-diagonal structure, which consists of a set of subproblems linked through a small number of state variables. Regarding the state variables in our formulation, we can see that only the capacity module number decisions $S$ and the capacity level change decisions $Y$ at node $n$ are dependent on node $a(n)$ from the previous stage. Then, $S_{a(n)}$ can be expressed in terms of $Y_{a(n)}$ via the linking constraints (\ref{eq:6}), enabling us to retain only the binary variables $Y$ as state variables. Therefore, we have the state variables $Y$, while the rest are local variables. We then reformulate the problem using dynamic programming. At the root node (i.e., $n = 1$), we have:

\begin{align}
    & z^{*} = \min \ 
    \Bigg(
        \sum_{j \in \boldsymbol{J}} 
        \sum_{l_1 \in \boldsymbol{L}} 
        \sum_{l_2 \in \boldsymbol{L}} 
        c_{j l_1 l_2 n}^F  
        Y_{j l_1 l_2 1} 
        \nonumber \\
    & \qquad
        + \sum_{(j,j^{'}) \in \boldsymbol{\mathit{JJ}}} 
        c_{j j^{'} n}^A  
        F_{j j^{'} 1} 
        + \sum_{(i,j) \in \boldsymbol{\mathit{IJ}}} 
        c_{i j n}^T  
        X_{i j 1} 
        \nonumber \\
    & \qquad
        + \sum_{j \in \boldsymbol{J}} 
        c_{n}^O  
        R_{j 1} 
        + \sum_{\nu \in \boldsymbol{\mathcal{C}(1)}} 
        p_{1 \nu} 
        Q_{\nu}(Y_1)
    \Bigg)
    \label{eq:18}
\end{align}

\vspace{-0.05in}
s.t., (\ref{eq:2})-(\ref{eq:17}) where $n = 1$,\\
\vspace{-0.1in}

\noindent where $z^{*}$ denotes the optimal value of the problem and $p_{n \nu}$ represents the conditional probability from node $n$ to its child node $\nu$. In addition, $Q_{\nu}(Y_1)$ denotes the optimal value function at child node $\nu$ given the current information from state variables $Y_1$. We point out that this new formulation is equivalent to the original problem (\ref{eq:1})-(\ref{eq:17}). Furthermore, we define the optimal value function $Q_n(\cdot)$ at each node $n \in \boldsymbol{\mathcal{T}} \backslash \{1\}$ as follows:

\begin{align}
    & Q_n\!\left(Y_{a(n)}\right) = \min \ 
    \Bigg(
        \sum_{j \in \boldsymbol{J}} 
        \sum_{l_1 \in \boldsymbol{L}} 
        \sum_{l_2 \in \boldsymbol{L}} 
        c_{j l_1 l_2 n}^F  
        Y_{j l_1 l_2 n}
        \nonumber \\
    & \qquad
        + \sum_{(j,j^{'}) \in \boldsymbol{\mathit{JJ}}} 
        c_{j j^{'} n}^A  
        F_{j j^{'} n}
        + \sum_{(i,j) \in \boldsymbol{\mathit{IJ}}} 
        c_{i j n}^T  
        X_{i j n}
        \nonumber \\
    & \qquad
        + \sum_{j \in \boldsymbol{J}} 
        c_{n}^O  
        R_{j n}
        + \sum_{\nu \in \boldsymbol{\mathcal{C}(n)}} 
        p_{n \nu} 
        Q_{\nu}(Y_n)
    \Bigg)
    \label{eq:19}
\end{align}

\vspace{-0.05in}
s.t., (\ref{eq:2})-(\ref{eq:17})\\
\vspace{-0.1in}

\noindent where the expected cost-to-go function at node $n$ can be denoted as $\phi_n(\cdot) := \sum_{\nu \in \boldsymbol{\mathcal{C}(n)}} p_{n \nu} Q_{\nu}(\cdot)$. 

The dynamic programming formulation above reveals the nested structure that we exploit in our SDDiP solution approach. Instead of attempting a direct and intractable backward induction, SDDiP iteratively builds approximations of the future costs. The algorithm alternates between a forward pass, where it proposes stage-feasible policies through sampled scenarios, and a backward pass, where it refines convex surrogates of the cost-to-go functions. These surrogates are gradually improved by generating backward cuts defined on the binary state variables $Y$, allowing the method to effectively tame the scenario-tree explosion while preserving the adaptivity at revision points.

\subsection{SDDiP Algorithm}
Prior to applying the SDDiP algorithm, it is essential to verify that the proposed formulation satisfies its underlying theoretical assumptions. In line with \cite{zou2019stochastic, zou2018multistage}, we examine the following key assumptions sequentially:

{
\renewcommand{\baselinestretch}{1.35}\selectfont
\begin{enumerate}[itemsep=0.8mm,parsep=0.05mm]
  \item {\textit{\textbf{Stage-wise independent scenario tree.}}} {We assume that the uncertain parameters, including demand $d_{in}$ and module throughput rate $k_{jn}$, follow a memory-less process. This ensures that the probability $p_n$ of any scenario node $n$ is independent of past events, satisfying the stage-wise independence.} 
  \item {\textit{\textbf{Binary state variables.}}} {As outlined in Section \ref{sec:dp}, the only state variable in our problem is the capacity level decision $Y$. According to constraints \eqref{eq:15}, $Y$ is strictly binary.}
  \item {\textit{\textbf{Markovian dependence.}}} {The formulation in Section \ref{sec:formulation} indicates that decisions at any node $n$ depend only on the state of the immediate parent node $a(n)$, as shown in state transition constraints \eqref{eq:5} and \eqref{eq:8}. }
  \item {\textit{\textbf{Complete continuous recourse.}}} {Our model guarantees feasibility through the inclusion of the continuous outsourcing variable $R$ in constraint \eqref{eq:4}. For any state $Y$ and set of integer decisions $\{F, S\}$, $R$ can be adjusted to meet the required demand, preventing infeasibility at the cost of outsourcing penalties.}
  \item {\textit{\textbf{Linearity and convexity.}}} {The formulation is a mixed-integer linear program. The objective function \eqref{eq:1} is a linear combination of the decision variables, and all constraints \eqref{eq:2}-\eqref{eq:10} are linear. This structure ensures that the problem at each stage is a nonempty, compact mixed-integer polyhedral set.}
\end{enumerate}
}

With these assumptions verified, we can then apply the SDDiP algorithm to solve the proposed problem. In the following subsections \ref{sec:fs}-\ref{sec:cg}, we develop a customized SDDiP for the proposed PAMSSP formulation. 

\subsubsection{Forward Pass}
\label{sec:fs}
Each iteration of SDDiP starts with scenario sampling. A subset of $M$ scenarios is randomly selected to represent possible realizations of uncertainty across stages. These scenarios form a set of paths from the root node to a leaf node in the scenario tree, capturing distinct outcomes for the uncertain parameters. Let $\boldsymbol{\Omega^k} = \{ \boldsymbol{\omega^k_{1}}, ..., \boldsymbol{\omega^k_{M}}\}$ denote the set of sampled scenarios at iteration $k$, where each $\boldsymbol{\omega^k_m}$ is a root-to-leaf path.

According to stage-wise independent assumption, we have $Q_n(\cdot) \equiv Q_t(\cdot),  \forall n \in \boldsymbol{\mathcal{T}_t}$. To formulate the SDDiP forward subproblem, we introduce auxiliary variables $Y^{P}$ as local copies of the state variables $Y_{a(n)}$, linking each node to its parent’s state. Moreover, these auxiliary variables are relaxed to be continuous. The forward subproblem for a given stage $t$ at iteration $k$ can be written as:

{
\begin{align}
    & \underline{Q}^k_t
    \!\left(
        \hat{Y}^k_{t-1},
        \underline{\phi}^k_t,
        \boldsymbol{\omega^{k}_{t m}}
    \right)
    = \nonumber \\
    & \min \ 
    \Bigg(
        \sum_{j \in \boldsymbol{J}} 
        \sum_{l_1 \in \boldsymbol{L}} 
        \sum_{l_2 \in \boldsymbol{L}} 
        c_{j l_1 l_2 n}^F  
        Y_{j l_1 l_2 n}
        \nonumber \\
    & \qquad
        + \sum_{(j,j^{'}) \in \boldsymbol{\mathit{JJ}}} 
        c_{j j^{'} n}^A  
        F_{j j^{'} n}
        + \sum_{(i,j) \in \boldsymbol{\mathit{IJ}}} 
        c_{i j n}^T  
        X_{i j n}
        \nonumber \\
    & \qquad
        + \sum_{j \in \boldsymbol{J}} 
        c_{n}^O  
        R_{j n}
        + \underline{\phi}^k_t(Y_t)
    \Bigg)
    \label{eq:20}
\end{align}

s.t., (\ref{eq:2})-(\ref{eq:17}) where $n \in \boldsymbol{\mathcal{T}_t} \cap \boldsymbol{\omega^{k}_{t m}}$ and $(\ref{eq:5}),(\ref{eq:8})$ are replaced by 

\vspace{0.1in}
\begin{align}
    & S_{jn}
    = \sum_{l_1 \in \boldsymbol{L}} 
      \sum_{l_2 \in \boldsymbol{L}} 
      u_{j l_2}  Y^{P}_{j l_1 l_2 n}  \nonumber \\
    & \qquad
      + \sum_{j^{'} \in \boldsymbol{J} \cup \{0\} \mid (j,j^{'}) \in \boldsymbol{\mathit{JJ}}} 
      \left(
          F_{j^{'}jn} - F_{jj^{'}n}
      \right), \nonumber \\
    & \qquad \qquad 
      \forall j \in \boldsymbol{J},\ 
      n \in \boldsymbol{\mathcal{T}_{t \backslash \{1\}}} 
      \cap \boldsymbol{\omega^{k}_{t m}}
      \label{eq:5.5} \\[0.4em]
    & \sum_{l_1 \in \boldsymbol{L}} 
      Y^{P}_{j, l_1, l, n}
    = \sum_{l_2 \in \boldsymbol{L}} 
      Y_{j, l, l_2, n}, \nonumber \\
    &  \qquad \qquad
      \forall j \in \boldsymbol{J},\ 
      l \in \boldsymbol{L},\ 
      n \in \boldsymbol{\mathcal{T}_{t \backslash \{1\}}} 
      \cap \boldsymbol{\omega^{k}_{t m}}
      \label{eq:8.5} 
\end{align}

\begin{align}
    & Y^{P}_{j l_1 l_2 n}
    = \hat{Y}^k_{j, l_1, l_2, t-1}, \nonumber \\
    &  \qquad \qquad
      \forall j \in \boldsymbol{J},\ 
      l_1, l_2 \in \boldsymbol{L},\ 
      n \in \boldsymbol{\mathcal{T}_t} 
      \cap \boldsymbol{\omega^{k}_{t m}}
      \label{eq:21} \\[0.4em]
    & Y^{P}_{j l_1 l_2 n} \in [0,1], \nonumber \\
    &  \qquad \qquad
      \forall j \in \boldsymbol{J},\ 
      l_1, l_2 \in \boldsymbol{L},\ 
      n \in \boldsymbol{\mathcal{T}_t} 
      \cap \boldsymbol{\omega^{k}_{t m}}
      \label{eq:22}
\end{align}
}

Although constraints (\ref{eq:21}) are redundant due to the introduction of copy variables, they play a critical role in supporting backward cut generation within the SDDiP framework \citep{zou2018multistage}. In addition, we note that $\underline{\phi}^k_t(Y_t)$ represents the lower-approximate expected cost-to-go value at stage $t$ given state information $Y_t$ and can be defined as the lower bound of a set of cutting planes:
\begin{align}
\underline{\phi}^k_t(Y_t) = \nonumber\\
\min \Bigg\{
    \theta_t :\ 
    & \theta_t \geq 
    \frac{1}{N_{t+1}} 
    \sum_{m=1}^{N_{t+1}} 
    \left[
        \alpha_{t+1, m}^{k^{'}} 
        + \left( \beta_{t+1, m}^{k^{'}} \right)^{\intercal} 
        Y_t
    \right] \nonumber \\
    & \text{and } \theta_t \geq L_t \ \text{for } \forall k^{'} < k
\Bigg\}
\label{eq:22.5}
\end{align}

\noindent where the $\theta_t$ is an auxiliary variable representing the estimated future cost, whose value is constrained by the collection of cutting planes generated in previous iterations. $L_{t}$ serves as a lower bound for the cost-to-go function to prevent unboundedness and is always set to zero in our problem. $N_{t}$ is the total number of realizations in stage $t$. $\alpha$ and $\beta$ are the corresponding cut coefficients attained from the backward pass.

For each sampled scenario $\boldsymbol{\omega^{k}_{m}} \in \boldsymbol{\Omega^{k}}$, the forward pass is solved from stage $1$ to stage $T$. Then, we can obtain a solution set $\{(\hat{Y}^k_{t m}, \hat{F}^k_{t m}, \hat{S}^k_{t m}, \hat{X}^k_{t m}, \hat{R}^k_{t m}), \forall t \in \boldsymbol{T} \}$ and the corresponding optimal objective value $z^k_{t m}$ for each stage $t$. Let $z^k_m := \sum_{t \in \boldsymbol{T}} \left( z^k_{t m} + \hat{\theta}^k_{tm} \right)$ be the modified objective value over $T$ stages. After finishing the forward pass for all sampled scenarios in $\boldsymbol{\Omega^{k}}$, the statistical upper bound $UB^k$ can be calculated as follows:

{
\renewcommand{\baselinestretch}{0.7}\selectfont
\begin{align}
&UB^k = \bar{z}^k + \mathit{z}_{\alpha/2} \cdot \frac{\sigma^k}{\sqrt{M}}, \text{ where } \bar{z}^k = \frac{1}{M} \sum_{m=1}^{M} z^k_m \\
& \qquad \qquad \text{ and } (\sigma^k)^2 = \frac{1}{M-1} \sum_{m=1}^{M} (z^k_m - \bar{z}^k)^2.
\end{align}
}

\subsubsection{Backward Pass with Cut Generation}
\label{sec:cg}
In the backward pass, feasible solutions from the forward pass are used to refine the value-function approximations from stage $T$ back to stage $1$ over the full scenario set $\boldsymbol{\Omega}$. At stage $T$ we set $\underline{\phi}^k_T(\cdot) \equiv 0$, and at each stage $t$ we solve relaxations of the forward subproblems at child nodes and use the resulting dual information to generate cuts, which update $\underline{\phi}^{k}_t(\cdot)$ to $\underline{\phi}^{k+1}_t(\cdot)$. The optimal value of the forward subproblem at stage $1$ is a valid lower bound, so we set $LB^k = \underline{Q}_{1}^{k}(v, \underline{\phi}^{k+1}_{1})$, where $v$ is the initial capacity level. To operationalize this backward update, we introduce the following cutting planes.
\vspace{0.05in}

\noindent \textbf{Benders Cut.} This type of cut can be easily obtained by solving the LP relaxation of the forward subproblem. Let $R_{\nu}^{k}(\hat{Y}_n^k, \underline{\phi}^{k+1}_{t_\nu}, \boldsymbol{\omega(\nu)})$ denote the LP relaxation of the forward subproblem at child node $\nu$, given the state variables $\hat{Y}_n^k$ obtained from the $k$-th iteration of the forward pass at the parent node $n$. Let $\gamma_{\nu}^{k}$ be the optimal value of $R_{\nu}^{k}(\hat{Y}_n^k, \underline{\phi}^{k+1}_{t_\nu}, \boldsymbol{\omega(\nu)})$. The Benders cut is then formulated as:

\begin{equation}
\theta_n \geq \sum_{\nu \in \boldsymbol{\mathit{C}(n)}} p_{n \nu} \left( \gamma_{\nu}^{k} + \pi_{\nu}^{k} {}^{{}^{\intercal}} (Y_n - \hat{Y}_n^k)\right),
\label{eq:23}
\end{equation}

\noindent where $\pi_{\nu}^k$ denotes the optimal dual value corresponding to linking constraints $Y^{P}_{\nu} = \hat{Y}^k_{n}$. Each time a cut is generated, the approximate expected cost-to-go function is updated. Accordingly, we use $\underline{\phi}^{k+1}_{t_\nu}$ instead of $\underline{\phi}^{k}_{t_\nu}$ in $R_{\nu}^k$ to reflect the most recent approximation. In addition, we remark that node $n$ always belongs to the sampling scenario set $\boldsymbol{\Omega^k}$, while the subproblems are solved for all of its child nodes $\boldsymbol{\mathit{C}(n)}$ to generate backward cuts.
\vspace{0.05in}

\noindent \textbf{Strengthened Benders Cut.} This cut is at least as tight as the classic Benders cut. Let $L^k_\nu(\hat{Y}^k_{n}, \underline{\phi}^{k+1}_{t_\nu}, \boldsymbol{\omega(\nu)})$ be the Lagrangian relaxation problem with penalizing the violation of linking constraints (\ref{eq:21}). Let $\pi_{\nu}^{k}$ denotes the optimal dual values associated with the linking constraints (\ref{eq:21}) obtained from solving $R_{\nu}^{k}(\hat{Y}_n^k, \underline{\phi}^{k+1}_{t_\nu}, \boldsymbol{\omega(\nu)})$. Then, the strengthened Benders cut takes the form:

\begin{equation}
\theta_n \geq \sum_{\nu \in \boldsymbol{\mathit{C}(n)}} p_{n \nu} \left( \eta_{\nu}^{k} + \pi_{\nu}^{k} {}^{{}^{\intercal}} Y_n \right),
\label{eq:24}
\end{equation}

\noindent where $\eta_{\nu}^{k}$ denotes the optimal objective value of the Lagrangian relaxation of the subproblem, evaluated using $\pi_{\nu}^{k}$ as fixed dual variables. 
\vspace{0.05in}

\noindent \textbf{Integer Optimality Cut.} Since all state variables in our problem are binary, the integer optimality cut can be formulated as follows:
\begin{align}
\theta_n \geq &\ (\bar{\gamma}_{n}^{k} - L_t) \Bigg( 
\sum_{j \in \boldsymbol{J}} \sum_{l_1 \in \boldsymbol{L}} \sum_{l_2 \in \boldsymbol{L}}\left(\hat{Y}_{j l_1 l_2 n}^k-1\right) Y_{j l_1 l_2 n} \nonumber\\
&+ \sum_{j \in \boldsymbol{J}} \sum_{l_1 \in \boldsymbol{L}} \sum_{l_2 \in \boldsymbol{L}}\left(Y_{j l_1 l_2 n}-1\right) \hat{Y}_{j l_1 l_2 n}^k \Bigg) + \bar{\gamma}_{n}^{k}, 
\end{align}

\noindent where $\bar{\gamma}_{n}^{k} = \sum_{\nu \in \boldsymbol{\mathit{C}(n)}} p_{n \nu} \gamma_{\nu}^{k}$ and $\gamma_{\nu}^{k}$ is the optimal value of forward subproblem $Q_{\nu}^{k}(\hat{Y}_n^k, \underline{\phi}^{k+1}_{t_\nu}, \boldsymbol{\omega(\nu)})$.

\vspace{0.05in}
\noindent \textbf{Lagrangian Cut.} This type of cut can be generated by solving the Lagrangian dual problem of the forward subproblem by dualizing the linking constraints (\ref{eq:21}). Let $D^k_\nu(\hat{Y}^k_{n}, \underline{\phi}^{k+1}_{t_\nu}, \boldsymbol{\omega(\nu)})$ be the Lagrangian dual problem. Let $\delta^k_{\nu}$ denote the optimal objective value of the Lagrangian dual problem, and $\lambda^k_{\nu}$ be the optimal Lagrangian multipliers corresponding to the linking constraints. Then, the Lagrangian cut can be written as follows:

\begin{equation}
\theta_n \geq \sum_{\nu \in \boldsymbol{\mathit{C}(n)}} p_{n \nu} \left( \delta_{\nu}^{k} + \lambda_{\nu}^{k} {}^{{}^{\intercal}} Y_n \right).
\end{equation}

The Lagrangian cuts are proven to be tight since the Lagrangian dual has zero duality gap when the state variables are binary \citep{zou2019stochastic}.

\subsection{SDDiP enhancements}
\subsubsection{Strengthened Cut Generation}
\label{sec:scg}

\cite{magnanti1981accelerating} proposed a method for generating Pareto-optimal cuts to improve the convergence of Benders decomposition. Instead of adding all potential cuts from the dual subproblem, this approach generates one of the strongest cuts in each iteration, which more efficiently prunes the master problem's feasible region and thereby accelerates the algorithm. A key element of this procedure is the use of a core point, representing a relative interior point of the master problem’s feasible region, which guides the generation of these high-quality cuts. After identifying a core point, the Magnanti-Wong secondary subproblem can be constructed and solved after solving the original dual subproblem to generate Pareto-optimal cuts. Building on this idea, \cite{papadakos2008practical} further introduced a more practical iterative scheme: at each iteration, the core point is updated by a convex combination of the previous core point $X^{\circ}$ and the current relaxed master solution $\hat{X}$, as $X^{\circ}_{new} \gets (1-\lambda)\cdot X^{\circ}_{old}+\lambda\cdot \hat{X}$, with $\lambda = 0.5$ based on empirical results. The author further stated that by updating the core point in this way, Pareto-optimal cuts could be directly derived from solving the independent Magnanti-Wong subproblems, which avoids the need to solve the original dual subproblem. This enhancement reduces computational complexity and addresses the numerical instability issue in solving the Magnanti-Wong secondary subproblem. To differentiate this approach from the classical Pareto-optimal cut generation, we refer to it as the independent Magnanti-Wong cuts.

Originally developed for Benders decomposition, these cut-generation methods are here adapted to the SDDiP framework for PAMSSP. The key modification concerns how core points are defined and updated across stages and nodes. Instead of a single master core point, SDDiP maintains core points for each stage–node pair so that they evolve along the scenario tree. In the forward pass, whenever a new solution is obtained at $(t,n)$, the algorithm either initializes the corresponding core point as the current solution if none exists, or updates it via $X^{\circ}_{t,n,\text{new}} \gets (1-\lambda)\cdot X^{\circ}_{t,n,\text{old}}+\lambda\cdot \hat{X}_{t,n}$ with $\lambda = 0.5$. These stage–node core points are then used to construct the Magnanti–Wong secondary and independent Magnanti–Wong subproblems that generate the strong cuts defined below.

\subsubsection*{Pareto-Optimal Cut.}
Let $R_{\nu}^{k}(\hat{Y}_n^k, \underline{\phi}^{k+1}_{t_\nu}, \boldsymbol{\omega(\nu)})$ be the LP relaxation of the forward subproblem. Then, using the core points $Y^{\circ k}_n$, we can modify $R_{\nu}^{k}$ so that its dual problem corresponds to the Magnanti-Wong secondary subproblem. We denote this new formulation as $M_{\nu}^{k}(\hat{Y}_n^k, Y^{\circ k}_n, \underline{\phi}^{k+1}_{t_\nu}, \boldsymbol{\omega(\nu)})$. Solving this modified problem yields the optimal dual multipliers $\alpha^{k}_{\nu}$ corresponding to the linking constraints (\ref{eq:21}), together with the optimal objective value $\rho_{\nu}^{k}$. The resulting Pareto-optimal cut is then given by:

\begin{equation}
\theta_n \geq \sum_{\nu \in \boldsymbol{\mathit{C}(n)}} p_{n \nu} \left( \rho_{\nu}^{k} + \alpha^{k}_{\nu} {}^{{}^{\intercal}} (Y_n - Y^{\circ k}_n) \right).
\label{eq:27}
\end{equation}

\noindent By construction, this cut leverages not only the current solution but also the historical information.

\subsubsection*{Independent Magnanti-Wong Cut.}
We further modify the subproblem $R_{\nu}^{k}$ such that its dual corresponds to the independent Magnanti-Wong subproblem. We denote it as $I_{\nu}^{k}(\hat{Y}_n^k, Y^{\circ k}_n, \underline{\phi}^{k+1}_{t_\nu}, \boldsymbol{\omega(\nu)})$. By solving $I_{\nu}^{k}$, we obtain the optimal dual multipliers $\beta^{k}
_{\nu}$ for linking constraints, and the optimal value $\zeta_{\nu}^{k}$. Then, the independent Magnanti-Wong cut is formulated by:

\begin{equation}
\theta_n \geq \sum_{\nu \in \boldsymbol{\mathit{C}(n)}} p_{n \nu} \left(\zeta_{\nu}^{k} + \beta^{k}_{\nu} {}^{{}^{\intercal}} (Y_n - Y^{\circ k}_n) \right).
\label{eq:28}
\end{equation}

\subsubsection*{Strengthened Pareto-Optimal Cut.}
To further accelerate convergence, we extend the Strengthened Benders cut of \cite{zou2019stochastic} to develop the strengthened Pareto-optimal and independent Magnanti-Wong cuts. The key idea is to integrate the dual information from the LP relaxation of the subproblem into a Lagrangian relaxation, where the linking constraints are relaxed and represented as penalty terms in the objective. The strengthened versions can produce tighter approximations of the cost-to-go function. Moreover, compared with pure Lagrangian cuts, the strengthened cuts avoid instability from oscillating subgradients and reduce the computational complexity. To implement this, we first solve the modified LP-relaxation subproblem and obtain the optimal dual multipliers of the linking constraints, just as in the Pareto-optimal or independent Magnanti-Wong cut generation. We then embed these multipliers into the Lagrangian relaxation problem to compute the optimal value and formulate the corresponding strengthened cut. By incorporating the optimal dual multipliers $\alpha_{\nu}^{k}$, obtained from solving $M_{\nu}^{k}$, into the Lagrangian relaxation problem $L^k_\nu(\hat{Y}^k_{n}, \underline{\phi}^{k+1}_{t_\nu}, \boldsymbol{\omega(\nu)})$, we derive its optimal value $\eta_{\nu}^k$. The strengthened Pareto-optimal cut is formulated as:

\begin{equation}
\theta_n \geq \sum_{\nu \in \boldsymbol{\mathit{C}(n)}} p_{n \nu} \left( \eta_{\nu}^{k} + \alpha^{k}_{\nu} {}^{{}^{\intercal}} Y_n \right).
\label{eq:29}
\end{equation}

\subsubsection*{Strengthened Independent Magnanti-Wong Cut.}
The optimal dual multipliers $\beta_{\nu}^{k}$ obtained from $I_{\nu}^{k}$ are incorporated into the Lagrangian relaxation $L^k_\nu$, from which we compute the optimal value $\eta_{\nu}^k$. The strengthened independent Magnanti-Wong cut can be written as:

\begin{equation}
\theta_n \geq \sum_{\nu \in \boldsymbol{\mathit{C}(n)}} p_{n \nu} \left(\eta_{\nu}^{k} + \beta^{k}_{\nu} {}^{{}^{\intercal}} Y_n \right).
\label{eq:30}
\end{equation}

\subsubsection*{{Validity for Strengthened Cuts}}

In our implementation of SDDiP, several enhanced cutting-plane families are used to approximate the value functions and play a central role in accelerating convergence. These cuts are obtained by evaluating a Lagrangian relaxation of the forward subproblems at the same dual multipliers that generate the classical version of cuts, but over a more restrictive feasible region. Intuitively, they preserve the dual slope while increasing the intercept, and therefore cannot produce a weaker lower bound.


\begin{proposition}[Dominance of Strengthened Pareto-Optimal Cuts]\label{prop:dominance}
Let $\hat{Y}_n^k$ and $Y^{\circ k}_n$ denote the state and core point selected at iteration $k$ during the forward pass. For each child node $\nu \in \mathcal{C}(n)$, let $\rho_\nu^k$ and $\alpha_\nu^k$ denote the optimal value and optimal dual multipliers of the Magnanti-Wong secondary problem $M_\nu^k(\hat{Y}_n^k, Y^{\circ k}_n)$ enforcing $Y^P_\nu = \hat{Y}_n^k$.
Let $\eta_\nu^k$ be the optimal value of the associated Lagrangian relaxation constructed using the same multipliers $\alpha_\nu^k$.
Then, for every feasible parent-state vector $Y_n$,

\begin{align}
 \sum_{\nu \in \boldsymbol{\mathit{C}(n)}} p_{n \nu} &\Big( \eta_{\nu}^{k} + (\alpha^{k}_{\nu})^{\intercal}  Y_n \Big) \notag\\
&\;\ge\;
\sum_{\nu \in \boldsymbol{\mathit{C}(n)}} p_{n \nu} \Big( \rho_{\nu}^{k} + (\alpha^{k}_{\nu})^{\intercal} (Y_n - Y^{\circ k}_n) \Big).
\end{align}

Hence, the strengthened Pareto-optimal cut is never weaker than classical Pareto-optimal cut.
\end{proposition}

A detailed proof is provided in Appendix~\ref{appendix:strengthened-benders}. A similar dominance proposition can also be formulated for strengthened independent Magnanti-Wong cuts. These dominance results not only guarantee the validity of the enhanced cuts but also help explain the faster lower-bound improvement observed in our computational experiments, especially when the LP relaxation exhibits a nonzero integrality gap.

\subsubsection{Alternating Cut Strategy}
\label{sec:acs}

Strengthened and Lagrangian cuts are tighter but more computationally expensive than LP-based cuts. To balance strength and effort, we design an alternating cut strategy in the SDDiP backward pass, inspired by \cite{angulo2016improving, arslan2021distribution}. For each state $(t,n,y^*)$, we maintain memory sets that record whether the state has been used to generate LP-based cuts (Benders-type, Pareto-optimal, independent Magnanti–Wong) or integer/Lagrangian cuts (integer optimality, Lagrangian, strengthened variants). When a state is first visited, LP-based cuts are generated; when revisited, integer/Lagrangian cuts are generated. A threshold parameter controls when memory sets are cleared so states can re-enter the alternation as new dual information becomes available.
A detailed description of the alternating cut strategy, including the memory sets $V_{LP}$ and $V$, the acceptance counter, the threshold $\zeta$, and the full pseudo-code of Algorithm~\ref{alg:opt_alg}, is provided in Appendix~\ref{appendix_algo_1}.

\subsubsection{Parallelization}
\label{sec:p}
In the SDDiP algorithm, the need to solve numerous subproblems at each stage in both the forward and backward passs presents a potential opportunity for acceleration through parallelization. We develop a general and scalable parallelization framework that exploits the independence of subproblem-solving tasks in both forward and backward passs, enabling efficient cut generation and aggregation. For each iteration, a master process samples scenario-tree nodes and partitions them evenly across processors to ensure balanced workloads. Stage-wise synchronization barriers coordinate progress, while cut coefficients generated in the backward pass are aggregated across processes to update the lower-approximation sets. Our implementation uses the Python \texttt{multiprocessing} library with shared dictionaries to synchronize data and preserve inter-process consistency.
The detailed parallelization procedure, including node partitioning rules, process synchronization, and pseudo-code of the parallel SDDiP implementation, is presented in Appendix~\ref{appendix_a}.

\section{Computational Study}
\label{sec:ne}
\subsection{Case Study}
\label{es}
We build a case study using a real data set from a three-year collaboration with a large U.S. modular construction company, which operates a hyperconnected network of modular assembly centers and mobile production containers. Modular construction relies on prefabricated building modules produced in factories and transported to sites, making total cost highly sensitive to shipping distance and the spatial deployment of capacity (Khan et al. 2024). Our partner company provides detailed monthly forecasts for upcoming modular housing projects over a three-year horizon, covering more than one hundred project locations across the United States. These forecasts are complemented by historical records on realized volumes and operational disruptions, which we use to calibrate the stochastic demand and supply processes in our scenario tree.

\begin{figure}[h!]
    \centering
    \includegraphics[width=\linewidth]{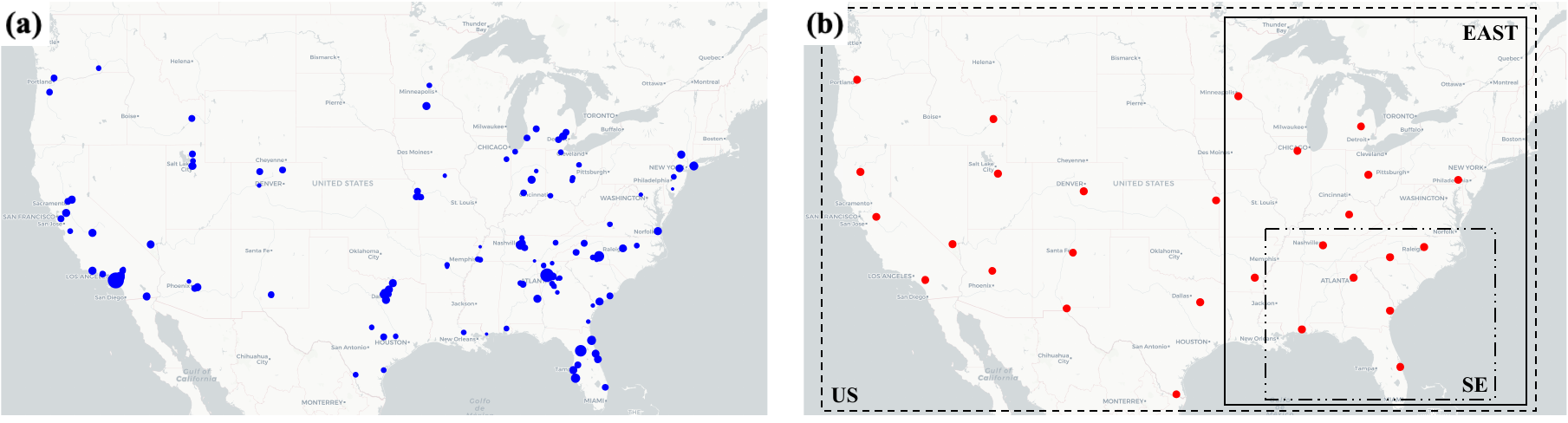}
    \caption{Illustration of (a) potential projects and (b) facility locations with three region instances.}
    \label{fig:network_info}
    
\end{figure} 

The data set includes the forecasted locations and monthly module demands of upcoming projects through 2026, together with a set of candidate assembly-center sites identified by the company as feasible leasing locations. Figure~\ref{fig:network_info} shows the project locations and candidate facilities, as well as the three regional instances; Table~\ref{tab:instance} summarizes their sizes. Across the three instances, the number of potential facilities ranges from 7 to 28 and the number of projects from 50 to 131, corresponding to a total demand between roughly 5{,}000 and 13{,}500 modules over the planning horizon. We consider:

{
\renewcommand{\baselinestretch}{0.9}\selectfont
\begin{itemize}
    \item Three nested regions (SE, EAST, US), capturing increasing network scale;
    \item Planning horizons of 3--12 months (monthly operational periods); 
    \item Two capacity-structure variants (2 and 3 capacity levels per facility).
\end{itemize}
}

Capacity expansion, reduction, and relocation costs are calibrated from company estimates and internal accounting data; all parameter values, demand and disruption scenarios, and the scenario-tree construction procedure are reported in Appendix~\ref{appendix_b}. The model is implemented in Python 3.9.12 with Gurobi 11.0.3 on a 24-core Intel Xeon Gold 6226 machine with 128 GB RAM. Subsequent subsections use this case to evaluate (i) the SDDiP enhancements, (ii) the value of partial adaptivity, and (iii) the value of stochastic, modular, and mobile capacity planning.

{
\vspace{-0.1in}
\renewcommand{\baselinestretch}{1.3}\selectfont
\begin{table}[h!]
\begin{center}
\caption{{Region instance characteristics.}}
\resizebox{0.4\linewidth}{!}{
\begin{tabular}{ccccc}
\toprule
Region & \#Facilities & \#Projects & \#Modules\\
\midrule
SE & 7 & 50 & 5,172 \\
EAST & 14 & 81 & 8,020 \\
US & 28 & 131 & 13,520 \\
\bottomrule
\end{tabular}
}
\label{tab:instance}
\end{center}
\vspace{-0.2in}
\end{table}
}

\subsection{Numerical Results}
\label{sec:nr}

To assess the computational performance of our proposed solution approach, we conduct two sets of numerical experiments. In the first experiment, we compare benchmark cutting planes from the SDDiP literature (i.e., integer optimality cuts (I), Lagrangian cuts (L), Benders cuts (B), and strengthened Benders cuts (SB)) against the new cuts proposed in this paper (i.e., Pareto-optimal cut (PT), the independent Magnanti-Wong cut (IM), the strengthened Pareto-optimal cut (SPT), and the strengthened independent Magnanti-Wong cut (SIM)), aiming to validate the effectiveness of the proposed cuts in Section \ref{sec:scg}. Cut performance is evaluated on the SE three-capacity-level instance under a 6-month planning horizon, with a two-hour runtime limit. Additional details on the experimental setup and results for the two-capacity-level instance are reported in Appendix~\ref{appendix_numerical_1}.

{
\renewcommand{\baselinestretch}{1.4}\selectfont
\begin{table}[h!]
\begin{center}
\caption{{Computational performance comparison across different cut combinations for a three-capacity-level instance.}}
\label{tab:exp1.1}
\resizebox{\linewidth}{!}{
\begin{tabular}{cccccccccccc}
\toprule
 \multicolumn{6}{c}{Classical SDDiP cutting planes} & \multicolumn{6}{c}{Proposed SDDiP cutting planes (Section \ref{sec:scg})} \\
\cmidrule(lr){1-6} \cmidrule(lr){7-12}
Cut type & \#Iter & \%Gap & \#Cuts & Runtime & Runtime/iter
& Cut type & \#Iter & \%Gap & \#Cuts & Runtime & Runtime/iter \\
\midrule
I      & 184 & 43.02 & 3,365 & 3,331 & 18.11
& PT + I   & 53 & 9.70 & 1,964 & 7,600 & 143.40 \\
L      & 58  & 11.40 & 1,069 & 7,578 & 130.66
& PT + L   & 33 & 0.19 & 1,226 & 6,889 & 208.77 \\
B + I  & 130 & 48.12 & 4,740 & 4,848 & 37.29
& IM + I   & 193 & 16.54 & 7,070 & 8,316 & 43.09 \\
B + L  & 48  & 18.07 & 1,784 & 7,594 & 158.21
& IM + L   & 63 & 0.19 & 2,310 & 7,641 & 121.28 \\
SB + I & 117 & 30.18 & 4,256 & 6,728 & 57.50
& SPT + I  & 18 & 0.19 &   664 & 1,456 & 80.91 \\
SB + L & 52  & 11.37 & 1,930 & 7,586 & 145.88
& SPT + L  & 21 & 0.19 &   778 & 4,032 & 192.01 \\
        &     &       &       &       &
& SIM + I  & 30 & 0.19 & 1,114 & 1,117 & 37.24 \\
        &     &       &       &       &
& SIM + L  & 23 & 0.19 &   856 & 2,752 & 119.66 \\
\bottomrule
\end{tabular}
}
\end{center}
\vspace*{-10pt}
\end{table}
}

Table \ref{tab:exp1.1} highlights the significant performance improvement of the newly proposed cutting planes compared to traditional benchmark cutting planes. The benchmark I and L cuts fail to achieve convergence, leaving large optimality gaps. Similarly, the introduction of B and SB cuts can accelerate convergence (i.e., reduce iteration numbers); however, neither is sufficient to close the optimality gap for larger instances. In contrast, the proposed cutting planes exhibit superior performance. Incorporating PT and IM cuts reduces remaining optimality gaps and improves computational efficiency, particularly combining with L cuts. Notably, the strengthened variants SPT and SIM achieve near-optimal solutions within only a few iterations, underscoring the effectiveness of the cutting planes. Among these, SPT and L cuts yield tighter gaps but generally require longer runtimes due to their higher per-iteration computational complexity. Overall, the combination of SIM and I cut combination provides the best trade-off between solution quality and runtime, achieving a 0.19\% optimality gap with the shortest runtime. Given its superior performance, we adopt the SIM + I cut configuration as the default strategy for all subsequent experiments. 

In the second experiment, we evaluate the efficiency and scalability of our enhanced SDDiP framework by comparing four algorithmic variants: (i) the classic sequential SDDiP (benchmark), (ii) a parallel SDDiP, (iii) a sequential SDDiP with the alternating cut strategy, and (iv) a parallel SDDiP with the alternating cut strategy. These algorithms are tested on the three-capacity-level instances of varying scale: three regions and three planning horizons with different scenario trees. Additional results for the two-capacity-level instances are provided in Appendix~\ref{appendix_numerical_2}. In this experiment, the runtime is limited to 24 hours, with parallelization carried out on 4 processors. 

{
\renewcommand{\baselinestretch}{1.3}\selectfont
\begin{table}[h!]
\begin{center}
\caption{{Computational performance of SDDiP algorithm variants for three-capacity-level instances.}}
\label{tab:exp1.2.2}
\resizebox{\linewidth}{!}{
\begin{tabular}{ccccccccccc}
\toprule
 & & & \multicolumn{4}{c}{Classic SDDiP} & \multicolumn{4}{c}{Parallel SDDiP} \\
\cmidrule(lr){4-7} \cmidrule(lr){8-11}
Region & \#Periods & \#Sce & \#Iter & \%Gap & Runtime & Runtime/Iter & \#Iter & \%Gap & Runtime & Runtime/Iter \\
\midrule
SE & 3 & 16 & 13 & 0.31 & 50 & 3.87 & 36 & 0.31 & 96 & 2.67 \\
& 6 & 1,024 & 23 & 0.19 & 478 & 20.79 & 36 & 0.19 & 639 & 17.74 \\
& 9 & 65,536 & 67 & 0.34 & 4,617 & 68.91 & 59 & 0.34 & 3,040 & 51.52 \\
EAST & 3 & 16 & 27 & 0.49 & 194 & 7.18 & 31 & 0.49 & 173 & 5.58 \\
& 6 & 1,024 & 50 & 0.15 & 3,016 & 60.33 & 48 & 0.15 & 2,311 & 48.16 \\
& 9 & 65,536 & 226 & 9.38 & 92,243 & 408.16 & 268 & 5.69 & 92,876 & 346.55 \\
US & 3 & 16 & 33 & 0.18 & 376 & 11.41 & 33 & 0.18 & 315 & 9.54 \\
& 6 & 1,024 & 233 & 2.72 & 91,945 & 394.61 & 163 & 2.68 & 52,687 & 241.68 \\
& 9 & 65,536 & 161 & 18.93 & 94,925 & 589.59 & 196 & 20.29 & 95,997 & 489.78 \\
\midrule
 & & & \multicolumn{4}{c}{SDDiP + alternating cut strategy} & \multicolumn{4}{c}{Parallel SDDiP + alternating cut strategy} \\
\cmidrule(lr){4-7} \cmidrule(lr){8-11}
Region & \#Periods & \#Sces & \#Iter & \%Gap & Runtime & Runtime/Iter & \#Iter & \%Gap & Runtime & Runtime/Iter \\
\midrule
SE & 3 & 16 & 14 & 0.31 & 37 & 2.65 & 36 & 0.30 & 52 & 1.44 \\
& 6 & 1,024 & 24 & 0.19 & 292 & 12.16 & 36 & 0.19 & 392 & 10.90 \\
& 9 & 65,536 & 52 & 0.42 & 1,700 & 32.70 & 61 & 0.58 & 2,204 & 36.13 \\
EAST & 3 & 16 & 25 & 0.49 & 106 & 4.23 & 33 & 0.49 & 116 & 3.51 \\
& 6 & 1,024 & 47 & 0.15 & 1,346 & 28.64 & 44 & 0.15 & 1,211 & 27.51 \\
& 9 & 65,536 & 367 & 3.46 & 92,374 & 251.70 & 334 & 7.00 & 91,504 & 273.96 \\
US & 3 & 16 & 40 & 0.18 & 283 & 7.08 & 53 & 0.18 & 345 & 6.51 \\
& 6 & 1,024 & 191 & 0.85 & 29,818 & 156.12 & 190 & 2.46 & 25,880 & 136.21 \\
& 9 & 65,536 & 225 & 14.57 & 95,996 & 426.65 & 237 & 17.53 & 96,414 & 406.81 \\
\bottomrule
\end{tabular}
}
\end{center}
\vspace*{-10pt}
\end{table}
}

Table~\ref{tab:exp1.2.2} reports the computational performance of different SDDiP variants for the three-capacity-level instances. We first compare the classic SDDiP and parallel SDDiP with their counterparts augmented by the alternating cut strategy. Across instances, incorporating alternating cuts leads to faster convergence and tighter optimality gaps, demonstrating that the alternating cut strategy consistently improves the effectiveness of the SDDiP algorithm under increased problem complexity. We then compare parallel SDDiP with classic sequential SDDiP. For most instances, parallelization reduces runtimes while maintaining similar optimality gaps, indicating that parallel SDDiP can effectively exploit computational resources. Finally, comparing the sequential SDDiP with alternating cuts against the parallel SDDiP with alternating cuts reveals that the sequential-alternating variant often achieves superior performance in three-capacity-level instances. In particular, for the EAST and US regions with longer planning horizons, the sequential-alternating approach yields tighter optimality gaps and comparable or even shorter runtimes. However, results for the two-capacity-level cases, summarized in Appendix~\ref{appendix_numerical_2}, exhibit a different pattern. In those smaller instances, the parallel SDDiP with alternating cuts outperforms the sequential variant in terms of runtime while achieving similar solution quality. Taken together, these findings indicate that while the alternating cut strategy robustly enhances SDDiP performance across problem scales, the effectiveness of parallelization is scale dependent. Large instances with higher complexity can diminish the benefit of parallelization due to increased memory requirements, communication overhead, and synchronization delays. Consequently, the choice between sequential and parallel variants with the alternating cut strategy should be guided by problem complexity. In the subsequent experiments, we therefore adopt the parallel SDDiP with the alternating cut strategy for the two-capacity-level cases, whereas the sequential-alternating approach is used for the three-capacity-level instances.

\subsection{The Value of Partially Adaptive Modeling}
\label{sec:pam}

This section assesses how partial adaptivity balances solution quality and computational complexity. We analyze several modeling configurations, starting with two baseline cases: (i) a model with a single revision point at the beginning of the planning horizon, which is equivalent to a TSSP, and (ii) a fully adaptive MSSP model with revisions in every period. We then gradually vary the number of revision points between these two extremes to change modeling configurations. To measure the corresponding adaptivity value, we adopt the value of partially adaptive multi-stage stochastic program (VPAMSP), a metric proposed by \cite{kayacik2025partially}:

\begin{equation}
VPAMSP = \frac{z^{\mathsf{TSSP}} - z^{\mathsf{PAMSSP}}(a)}{z^{\mathsf{TSSP}} - z^{\mathsf{MSSP}}} \times 100\%.
\end{equation}    

\noindent where $z$ is the objective value and $a$ is the number of revision points. This metric evaluates the proportion of the adaptivity value captured by the partially adaptive model PAMSSP compared to TSSP, normalized by the value gap between TSSP and fully adaptive MSSP. 

{
\renewcommand{\baselinestretch}{1.3}\selectfont
\begin{table}[h!]
\begin{center}
\caption{{Solving PAMSSP under various numbers of revision points.}}
\label{tab:exp2}
\resizebox{\linewidth}{!}{
\begin{tabular}{ccccccccccc}
\toprule
& & & \multicolumn{4}{c}{2 capacity levels} & \multicolumn{4}{c}{3 capacity levels} \\
\cmidrule(lr){4-7} \cmidrule(lr){8-11}
\#Periods & \#Rev & Revision points & \#Iter & \%Gap & Runtime & \%VPAMSP & \#Iter & \%Gap & Runtime & \%VPAMSP\\
\midrule
3 & 1 & 1       & 27  & 0.42 & 57 & 0.00 & 29  & 0.37 & 115 & 0.00 \\
  & 2 & 1,3     & 26  & 0.45 & 51 & 44.02 & 29  & 0.49 & 127 & 66.65 \\
  & 3 & 1,2,3   & 248 & 0.49 & 762 & 100.00 & 50  & 0.53 & 238 & 100.00 \\
\midrule
6 & 1 & 1       & 62  & 0.11 & 687 & 0.00 & 74  & 0.12 & 1,653 & 0.00 \\
  & 2 & 1,4     & 70  & 0.12 & 944 & 19.83 & 87  & 0.25 & 2,782 & 26.88 \\
  & 4 & 1,2,4,6 & 69  & 0.13 & 939 & 66.16 & 199 & 0.19 & 10,625 & 77.34 \\
  & 6 & 1,2,3,4,5,6 & 998 & 11.36 & 89,088 & 100.00 & 230 & 4.68 & 25,375 & 100.00 \\
\midrule
9 & 1 & 1       & 140 & 0.38 & 7,127 & 0.00 & 397 & 0.26 & 62,811 & 0.00 \\
  & 2 & 1,5     & 154 & 0.07 & 9,161 & 18.27 & 337 & 2.02 & 94,372 & 35.33 \\
  & 4 & 1,3,5,7 & 147 & 0.07 & 8,454 & 45.82 & 319 & 5.53 & 91,597 & 63.51 \\
  & 6 & 1,2,4,6,8,9 & 230 & 0.92 & 18,070 & 71.48 & 354 & 5.91 & 91,575 & 86.39 \\
  & 9 & 1,2,3,4,5,6,7,8,9 & 500 & 9.72 & 83,608 & 100.00 & 251 &  13.40 & 90,524 & 100.00 \\
\bottomrule
\end{tabular}
}
\end{center}
\vspace*{-8pt}
\end{table}
}

Table \ref{tab:exp2} presents the results for the EAST region instance with two and three capacity levels under varying planning horizons. As expected, TSSP models (i.e., \#Rev$ = 1$) solve fastest but fail to capture adaptivity. Increasing the number of revision points consistently improves solution quality, as reflected by higher VPAMSP, but at the cost of longer runtimes. This is more evident in MSSP models (\#Rev$ = T$), which achieve maximum adaptivity but often incur prohibitive runtimes or fail to converge within the 24-hour limit.

Another important observation is that a relatively small number of revision points enables the partially adaptive models to recover a large fraction of VPAMSP while keeping computational times within practical limits. For example, in the 6-period, 3-capacity-level case, setting revision points at \{1,2,4,6\} achieves 77.34\% VPAMSP while requiring less than half the runtime of the unsolved fully adaptive model. In addition, the results indicate that the marginal benefit of additional revision points decreases as more points are introduced. Moving from two to four revision points yields substantial gains of 50.31\% VPAMSP in the 6-period, 3-capacity-level case, whereas the incremental benefit of increasing from four revision points to full adaptivity is much smaller (i.e., 22.66\% VPAMSP) relative to the additional computational effort. This finding highlights the practical importance of partially adaptive modeling: a small number of well-placed revision points can capture most of the benefits of adaptivity with significantly lower computational effort.

\subsection{The Value of Dynamic Stochastic Modeling}
\label{sec:vss}

This section investigates the advantages of explicitly accounting for uncertainty through dynamic stochastic modeling, compared to its deterministic counterpart. We adopt the value of the stochastic solution metrics, including expected value (EV), expectation of expected value at stage $T$ (EEV$_T$), recourse problem value (RP), and value of the stochastic solution at stage $T$ (VSS$_T$) (see Appendix \ref{appendix_exp3}). The experiment in Table \ref{tab:exp3} is conducted across varying uncertainty levels by adjusting the demand standard deviation $\sigma$ and disruption rate $\lambda$ in the scenario tree generation.

\begin{table}[h!]
\begin{center}
\caption{{Multi-stage stochastic programming metrics under different uncertainty parameters.}}
\label{tab:exp3}
\resizebox{\linewidth}{!}{
\begin{tabular}{cccccccccc}
\toprule
& & \multicolumn{4}{c}{3 planning periods} & \multicolumn{4}{c}{6 planning periods} \\
\cmidrule(lr){3-6} \cmidrule(lr){7-10}
\quad $\sigma$ \quad\ & \quad $\lambda$ \quad\
& EV & EEV$_T$ & RP & VSS$_T$ 
& EV & EEV$_T$ & RP & VSS$_T$ \\
\midrule
0.5 & 0.5 & 3,020,376 & 4,246,648 & 4,233,056 & 13,592 
          & 4,545,280 & 5,319,681 & 5,208,949 & 110,732 \\
0.5 & 1.0 & 3,020,376 & 4,503,452 & 4,235,279 & 268,173 
          & 4,545,280 & 6,383,246 & 5,970,628 & 412,618 \\
0.5 & 1.5 & 3,020,376 & 5,694,260 & 5,203,824 & 490,436 
          & 4,545,280 & 8,593,500 & 8,192,850 & 400,650 \\
0.5 & 2.0 & 3,020,376 & 7,127,225 & 6,045,162 & 1,082,063 
          & 4,545,280 & 10,971,657 & 9,166,326 & 1,805,331 \\
1.0 & 0.5 & 3,085,480 & 4,331,740 & 4,266,701 & 65,039 
          & 4,591,981 & 5,629,275 & 5,252,453 & 376,822 \\
1.5 & 0.5 & 3,099,833 & 4,349,238 & 4,274,307 & 74,931 
          & 4,602,456 & 5,645,005 & 5,261,922 & 383,083 \\
2.0 & 0.5 & 3,104,544 & 4,355,814 & 4,273,938 & 81,876 
          & 4,592,470 & 5,647,013 & 5,261,696 & 385,317 \\
\bottomrule
\end{tabular}
}
\end{center}
\vspace*{-10pt}
\end{table}

Across all instances in Table~\ref{tab:exp3}, we observe that re-evaluating the EV decisions under under stochastic realizations yields high EEV$_T$ values. This indicates that deterministic solutions are overly optimistic and fail to hedge against uncertainty. In contrast, the RP values obtained from PAMSSP are consistently lower than EEV$_T$, confirming that explicitly modeling uncertainty enables better risk mitigation. Consequently, VSS$_T$ is always positive and measures the economic value of stochastic solutions. In addition, the value of stochastic modeling increases as the planning horizon extends from 3 to 6 periods. Longer horizons accumulate greater uncertainty, creating more opportunities for stochastic formulations to outperform deterministic policies. Accordingly, dynamic adaptivity delivers increasing benefits in mitigating uncertainty-induced losses.

Furthermore, the results show that the advantages of stochastic modeling rise with uncertainty intensity. For supply uncertainty, in the 3-period case with $\sigma=0.5$, VSS$_T$ increases from $\$13,592$ at $\lambda=0.5$ to over $\$1,082,063$ at $\lambda=2.0$, where the stochastic model captures more value. Similarly, higher demand volatility amplifies the value of stochastic recourse. Overall, results show that ignoring uncertainty leads to suboptimal decisions, and the proposed stochastic model delivers greater cost savings and robustness under demand and supply uncertainty.

\subsection{The Value of Modular and Mobile Facilities}
\label{sec:vmm}

This final experiment conducts a sensitivity analysis to quantify the economic benefits of modularity and mobility. We adopt the value of module mobility (VMM) \citep{allman2020dynamic} to measure the total rewards of adopting modular and mobile designs compared to the entirely static facility network (i.e., non-modular, non-mobile baseline):
\vspace{-0.01in}
\begin{equation}
VMM = z^*_{\text{Non-mod\&Non-mob}} - z^*_{\text{Mod\&Mob}},
\end{equation}

\noindent where $z^*$ denotes the optimal value of the proposed model under the corresponding setting.
We further propose two new metrics, named value of modularity (VMoD) and value of mobility (VMoB), to evaluate the individual added value of each feature. Among these, VMoD measures the cost savings from being able to adjust a facility's capacity in discrete, modular units, compared to a baseline where facility capacities are fixed: 

\begin{equation}
VMoD = z^*_{\text{Non-mod\&Non-mob}} - z^*_{\text{Mod\&Non-mob}}.
\end{equation}

\noindent While VMoB quantifies the advantage of having the option to relocate capacity modules between modular facilities, compared to a system limited to only adding or removing capacity units:

\begin{equation}
VMoB = z^*_{\text{Mod\&Non-mob}} - z^*_{\text{Mod\&Mob}}.
\end{equation}

\begin{figure}[h!]
    \centering
    \includegraphics[width=\linewidth]{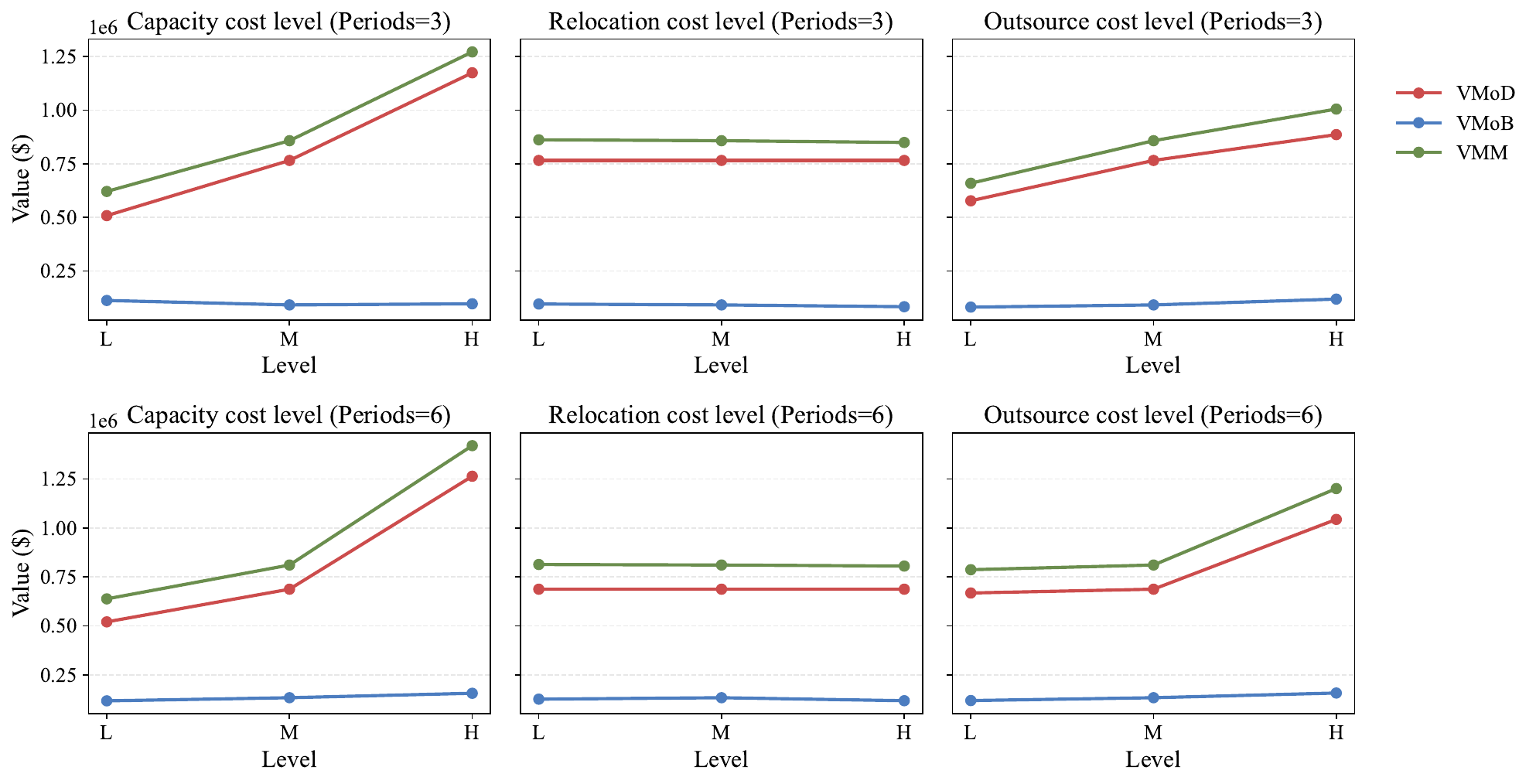}
    \caption{{Cost parameter sensitivity analysis under different modularity and mobility options.}}
    \label{fig:exp4}
    \vspace{-8pt}
\end{figure}

We evaluate the above metrics under various cost parameters: (i) facility capacity costs, (ii) module relocation costs, and (iii) outsourcing penalty costs. In Figure \ref{fig:exp4}, cost levels are classified as low (L), medium (M), and high (H), where M denotes the baseline, and L and H are set to half and double the baseline values. These variations in cost parameters reflect different regional configurations. For example, in economically developed and densely populated regions, capacity costs tend to be higher due to expensive labor costs and facility rental prices. While rural areas may face higher relocation costs, reflecting longer transport distances and underdeveloped infrastructure.

The full numerical results are reported in Appendix \ref{appendix_exp4}, while Figure \ref{fig:exp4} visualizes the metrics. As shown, introducing modularity alone yields substantial improvements over the non-modular baseline. In the 3-period case with medium cost parameters, VMoD equals \$765,878, indicating a cost reduction of over 15\%. These gains become more pronounced under high-outsourcing cost and high-capacity cost scenarios, where VMoD exceeds \$1.04 million and \$1.26 million, respectively, in the 6-month planning horizon, underscoring the importance of modular capacity adjustment. Mobility further enhances system performance, with VMoB ranging from approximately \$80,000 to \$160,000, indicating consistent additional savings by allowing the relocation of capacity modules across facilities. However, relative to modularity, the marginal gain of mobility is limited, accounting for roughly 2\% of the baseline cost. This reflects the trade-off inherent in relocation decisions: relocating an existing module is advantageous only when facilities are geographically close and demand shifts are well aligned in time; otherwise, returning unused modules and renting new capacity is more cost-effective. As a result, mobility provides flexibility in certain cases, while modularity delivers broader and more robust system-wide cost reductions.

\section{Conclusion}
\label{sec:con}
This paper addresses the Dynamic Stochastic Modular and Mobile Capacity Planning (DSMMCP) problem, a novel tactical optimization challenge emerging from the shift toward more distributed and hyperconnected supply chain networks. The DSMMCP problem leverages the short-term leasing of modular facilities and the dynamic relocation of mobile capacity modules, enabling companies to better respond to rapidly changing market conditions such as demand fluctuations and supply disruptions. We formulate DSMMCP as a partially adaptive multi-stage stochastic program (PAMSSP), balancing decision adaptability and computational complexity by pre-committing the timing of capacity adjustments while letting their magnitudes respond to realized information. To solve this challenging problem, we design an enhanced stochastic dual dynamic integer programming (SDDiP) algorithm with strengthened cut generation, an alternating cut strategy, and a parallelization framework tailored for SDDiP. To the best of our knowledge, this work represents the first application of SDDiP enhanced with the alternating cut strategy to address PAMSSPs. \\
We then conduct a case study based on a three-year collaboration with a large U.S. modular construction company, where the experimental instances reflect real project forecasts and potential facility sites. The computational results demonstrate the effectiveness of the proposed modeling and solution framework. The enhanced cut strategies significantly reduce runtime and deliver high-quality solutions compared to classic cutting planes. Among these, the combination of strengthened independent Magnanti-Wong cut and Lagrangian cut (SIM+L) performs the best, providing near-optimal solutions within practical runtimes, and the alternating cut strategy further accelerates convergence, with parallelization particularly beneficial for moderate-size instances. The experiments on partially adaptive modeling show that a PAMSSP model with a few well-placed revision points captures most of the value of full adaptivity, while the dynamic stochastic formulation yields a value of the stochastic solution (VSS) above \$0.4 million in the 6-period case and increasing under higher uncertainty. Sensitivity analysis further highlights the economic benefits of modularity and mobility in the 6-period high-capacity-cost case, where modularity alone yields savings above \$1.2 million, mobility adds about \$150{,}000, and the combined value of modular mobility (VMM) exceeds \$1.4 million.\\
There are several research avenues for future research. First, multi-objective extensions of DSMMCP could jointly address cost, service, resilience, and environmental performance, for example by internalizing carbon emissions or energy use in the objective. Second, combining partially adaptive modeling with rolling-horizon schemes would allow revision points to be updated dynamically as new information becomes available. Third, applications in disaster response and emergency supply chains, where rapid deployment and relocation of temporary facilities are critical, would further illustrate the potential of DSMMCP to enhance societal resilience under crisis conditions.


\vspace{-0.05in}
\bibliographystyle{elsarticle-harv} 
\bibliography{scibib}

@techreport{abb2024,
  author    = {ABB},
  year      = {2024},
  title     = {Ninety-one percent of industrial businesses hit by resource scarcity highlighting circularity needs},
  institution = {Asea Brown Boveri Motion},
}

@article{faugere2022dynamic,
  author    = {Faug{\`e}re, L. and Klibi, W. and White~III, C. and Montreuil, B.},
  year      = {2022},
  title     = {Dynamic pooled capacity deployment for urban parcel logistics},
  journal   = {European Journal of Operational Research},
  volume    = {303},
  pages     = {650--667}
}

@techreport{Mibox2024,
  author      = {MI-BOX},
  title       = {Mobile Storage Containers},
  institution     = {MI-BOX Mobile Storage},
  year        = {2024},
  type        = {Technical Report},
}

@article{kim2021hyperconnected,
  author = {Kim, N. and Montreuil, B. and Klibi, W. and Kholgade, N.},
  year = {2021},
  title = {Hyperconnected urban fulfillment and delivery},
  journal = {Transportation Research Part E: Logistics and Transportation Review},
  volume = {145},
  pages = {102104}
}

@article{kulkarni2022resilient,
  author = {Kulkarni, O. and Dahan, M. and Montreuil, B.},
  year = {2022},
  title = {Resilient hyperconnected parcel delivery network design under disruption risks},
  journal = {International Journal of Production Economics},
  volume = {251},
  pages = {108499}
}

@inproceedings{liu2023logistics,
  title={Logistics Hub Capacity Deployment in Hyperconnected Transportation Network Under Uncertainty},
  author={Liu, Xiaoyue and Li, Jingze and Montreuil, Benoit},
  booktitle = {IISE Annual Conference Proceedings},
  publisher = {Institute of Industrial and Systems Engineers (IISE)},
  year={2023}
}

@article{liu2025dynamic,
  title={Dynamic hub capacity planning in hyperconnected relay transportation networks under uncertainty},
  author={Liu, Xiaoyue and Li, Jingze and Dahan, Mathieu and Montreuil, Benoit},
  journal={Transportation Research Part E: Logistics and Transportation Review},
  volume={194},
  pages={103940},
  year={2025},
  publisher={Elsevier}
}

@article{gill2008identifying,
  title={Identify potential bottlenecks through activity underutilization costs},
  author={Gill, A},
  journal={International Journal of Simulation Modelling (IJSIMM)},
  volume={7},
  number={4},
  year={2008}
}

@article{magnanti1981accelerating,
  title={Accelerating Benders decomposition: Algorithmic enhancement and model selection criteria},
  author={Magnanti, Thomas L and Wong, Richard T},
  journal={Operations research},
  volume={29},
  number={3},
  pages={464--484},
  year={1981},
  publisher={INFORMS}
}

@article{papadakos2008practical,
  title={Practical enhancements to the Magnanti--Wong method},
  author={Papadakos, Nikolaos},
  journal={Operations Research Letters},
  volume={36},
  number={4},
  pages={444--449},
  year={2008},
  publisher={Elsevier}
}

@techreport{Symonds2019,
  author      = {OnTrac},
  title       = {Self-driving parcel lockers and delivery vehicles to be trialled in Northern California},
  institution     = {Parcel and Postal Technology International},
  year        = {2019},
  type        = {Technical Report},
}

@techreport{DefenceRedefined2023,
  author      = {Rheinmetall},
  year        = {2023},
  title       = {Presentation of the Mobile Smart Factory for Battle Damage Repair},
  institution = {Defence Redefined},
  type        = {Technical Report},
}

@article{montreuil2011toward,
  author = {Montreuil, B.},
  year = {2011},
  title = {Toward a Physical Internet: meeting the global logistics sustainability grand challenge},
  journal = {Logistics Research},
  volume = {3},
  pages = {71--87}
}

@book{montreuil2013foundations,
  title={Physical internet foundations},
  author={Montreuil, Benoit and Meller, Russell D and Ballot, Eric},
  year={2013},
  publisher={Springer}
}

@book{farahani2009facility,
  title={Facility location: concepts, models, algorithms and case studies},
  author={Farahani, Reza Zanjirani and Hekmatfar, Masoud},
  year={2009},
  publisher={Springer Science \& Business Media}
}

@article{zou2019stochastic,
  author = {Zou, J. and Ahmed, S. and Sun, X.},
  year = {2019},
  title = {Stochastic dual dynamic integer programming},
  journal = {Mathematical Programming},
  volume = {175},
  pages = {461--502}
}

@article{birge1985decomposition,
  title={Decomposition and partitioning methods for multistage stochastic linear programs},
  author={Birge, John R},
  journal={Operations research},
  volume={33},
  number={5},
  pages={989--1007},
  year={1985},
  publisher={INFORMS}
}

@book{birge1997introduction,
  title={Introduction to stochastic programming},
  author={Birge, John R and Louveaux, Francois},
  year={1997},
  publisher={Springer}
}

@book{shapiro2021lectures,
  title={Lectures on stochastic programming: modeling and theory},
  author={Shapiro, Alexander and Dentcheva, Darinka and Ruszczynski, Andrzej},
  year={2021},
  publisher={SIAM}
}

@article{lara2020electric,
  title={Electric power infrastructure planning under uncertainty: stochastic dual dynamic integer programming (SDDiP) and parallelization scheme},
  author={Lara, Cristiana L and Siirola, John D and Grossmann, Ignacio E},
  journal={Optimization and Engineering},
  volume={21},
  pages={1243--1281},
  year={2020},
  publisher={Springer}
}

@article{pereira1991multi,
  title={Multi-stage stochastic optimization applied to energy planning},
  author={Pereira, Mario VF and Pinto, Leontina MVG},
  journal={Mathematical programming},
  volume={52},
  pages={359--375},
  year={1991},
  publisher={Springer}
}

@article{quezada2022combining,
  title={Combining polyhedral approaches and stochastic dual dynamic integer programming for solving the uncapacitated lot-sizing problem under uncertainty},
  author={Quezada, Franco and Gicquel, C{\'e}line and Kedad-Sidhoum, Safia},
  journal={INFORMS Journal on Computing},
  volume={34},
  number={2},
  pages={1024--1041},
  year={2022},
  publisher={INFORMS}
}

@article{angulo2016improving,
  title={Improving the integer L-shaped method},
  author={Angulo, Gustavo and Ahmed, Shabbir and Dey, Santanu S},
  journal={INFORMS Journal on Computing},
  volume={28},
  number={3},
  pages={483--499},
  year={2016},
  publisher={INFORMS}
}

@article{zou2018multistage,
  title={Multistage stochastic unit commitment using stochastic dual dynamic integer programming},
  author={Zou, Jikai and Ahmed, Shabbir and Sun, Xu Andy},
  journal={IEEE Transactions on Power Systems},
  volume={34},
  number={3},
  pages={1814--1823},
  year={2018},
  publisher={IEEE}
}

@article{jena2015dynamic,
  title={Dynamic facility location with generalized modular capacities},
  author={Jena, Sanjay Dominik and Cordeau, Jean-Fran{\c{c}}ois and Gendron, Bernard},
  journal={Transportation Science},
  volume={49},
  number={3},
  pages={484--499},
  year={2015},
  publisher={INFORMS}
}

@article{jena2015modeling,
  title={Modeling and solving a logging camp location problem},
  author={Jena, Sanjay Dominik and Cordeau, Jean-Fran{\c{c}}ois and Gendron, Bernard},
  journal={Annals of Operations Research},
  volume={232},
  pages={151--177},
  year={2015},
  publisher={Springer}
}

@inproceedings{marcotte2016introducing,
  title        = {Introducing the Concept of Hyperconnected Mobile Production},
  author       = {Marcotte, Suzanne and Montreuil, Benoit},
  booktitle    = {14th IMHRC Proceedings (Karlsruhe, Germany – 2016)},
  year         = {2016},
  publisher    = {Progress in Material Handling Research},
}

@inproceedings{marcotte2015modeling,
  title={Modeling of physical internet enabled interconnected modular production},
  author={Marcotte, S and Montreuil, B and Coelho, LC},
  booktitle={Proceedings of 2nd International Physical Internet Conference (IPIC)},
  year={2015}
}

@article{jena2017lagrangian,
  title={Lagrangian heuristics for large-scale dynamic facility location with generalized modular capacities},
  author={Jena, Sanjay Dominik and Cordeau, Jean-Fran{\c{c}}ois and Gendron, Bernard},
  journal={INFORMS Journal on Computing},
  volume={29},
  number={3},
  pages={388--404},
  year={2017},
  publisher={INFORMS}
}

@article{becker2019value,
  title={Value of modular production concepts in future chemical industry production networks},
  author={Becker, Tristan and Lier, Stefan and Werners, Brigitte},
  journal={European Journal of Operational Research},
  volume={276},
  number={3},
  pages={957--970},
  year={2019},
  publisher={Elsevier}
}

@article{hong2020optimal,
  title={Optimal planning and modular infrastructure dynamic allocation for shale gas production},
  author={Hong, Bingyuan and Li, Xiaoping and Song, Shangfei and Chen, Shilin and Zhao, Changlong and Gong, Jing},
  journal={Applied Energy},
  volume={261},
  pages={114439},
  year={2020},
  publisher={Elsevier}
}

@article{correia2021integrated,
  title={Integrated facility location and capacity planning under uncertainty},
  author={Correia, Isabel and Melo, Teresa},
  journal={Computational and applied mathematics},
  volume={40},
  number={5},
  pages={175},
  year={2021},
  publisher={Springer}
}

@article{allman2020dynamic,
  title={Dynamic location of modular manufacturing facilities with relocation of individual modules},
  author={Allman, Andrew and Zhang, Qi},
  journal={European Journal of Operational Research},
  volume={286},
  number={2},
  pages={494--507},
  year={2020},
  publisher={Elsevier}
}

@book{vstadlerova2024solving,
  title={Solving Multi-stage Stochastic Facility Location Problems with Modular Capacity Adjustments},
  author={{\v{S}}t{\'a}dlerov{\'a}, {\v{S}}{\'a}rka and Sch{\"u}tz, Peter and Jena, Sanjay Dominik},
  year={2024},
  publisher={Bureau de Montreal, Universit{\'e} de Montreal}
}

@article{kayser2023relocatable,
  title={Relocatable modular capacities in risk aware strategic supply network planning under demand uncertainty},
  author={Kayser, Ariane and Sahling, Florian},
  journal={Schmalenbach Journal of Business Research},
  volume={75},
  number={1},
  pages={1--35},
  year={2023},
  publisher={Springer}
}

@article{correia2024matheuristic,
  title={A matheuristic for a multi-period three-echelon network design problem with temporary capacity acquisition},
  author={Correia, Isabel and Melo, Teresa},
  journal={Computers \& Industrial Engineering},
  volume={192},
  pages={110244},
  year={2024},
  publisher={Elsevier}
}

@article{ge2023multistage,
  title={Multistage stochastic programming for the closed-loop supply chain planning with mobile modules under uncertainty},
  author={Ge, Congqin and Zhang, Lifeng and Yang, Wenhui and Yuan, Zhihong},
  journal={AIChE Journal},
  volume={69},
  number={9},
  pages={e18156},
  year={2023},
  publisher={Wiley Online Library}
}

@article{kayacik2025partially,
  title={Partially adaptive multistage stochastic programming},
  author={Kayac{\i}k, Sezen Ece and Basciftci, Beste and Schrotenboer, Albert H and Ursavas, Evrim},
  journal={European Journal of Operational Research},
  volume={321},
  number={1},
  pages={192--207},
  year={2025},
  publisher={Elsevier}
}

@article{niaz2022revolutionizing,
  title={Revolutionizing inventory planning: Harnessing digital supply data through digitization to optimize storage efficiency pre-and post-pandemic},
  author={Niaz, Moazam},
  journal={BULLET: Jurnal Multidisiplin Ilmu},
  volume={1},
  number={03},
  pages={592273},
  year={2022}
}

@article{asawawibul2025influence,
  title={The influence of cost on customer satisfaction in e-commerce logistics: Mediating roles of service quality, technology usage, transportation time, and production condition},
  author={Asawawibul, Sasiprapha and Na-Nan, Khahan and Pinkajay, Kaptun and Jaturat, Nutt and Kittichotsatsawat, Yotsaphat and Hu, Bowei},
  journal={Journal of Open Innovation: Technology, Market, and Complexity},
  volume={11},
  number={1},
  pages={100482},
  year={2025},
  publisher={Elsevier}
}

@article{urban2025commerce,
  title={E-commerce enterprise flexibility leading to better customer perception},
  author={Urban, Wieslaw and Buraczy{\'n}ska, Barbara},
  journal={Journal of Retailing and Consumer Services},
  volume={85},
  pages={104267},
  year={2025},
  publisher={Elsevier}
}

@article{alarcon2022modular,
  title={Modular and mobile facility location problems: A systematic review},
  author={Alarcon-Gerbier, Eduardo and Buscher, Udo},
  journal={Computers \& Industrial Engineering},
  volume={173},
  pages={108734},
  year={2022},
  publisher={Elsevier}
}

@article{zou2018partially,
  title={Partially adaptive stochastic optimization for electric power generation expansion planning},
  author={Zou, Jikai and Ahmed, Shabbir and Sun, Xu Andy},
  journal={INFORMS journal on computing},
  volume={30},
  number={2},
  pages={388--401},
  year={2018},
  publisher={INFORMS}
}

@article{basciftci2024adaptive,
  title={Adaptive Two-Stage Stochastic Programming with an Analysis on Capacity Expansion Planning Problem},
  author={Basciftci, Beste and Ahmed, Shabbir and Gebraeel, Nagi},
  journal={Manufacturing \& Service Operations Management},
  year={2024},
  publisher={INFORMS}
}

@article{arslan2021distribution,
  title={Distribution network deployment for omnichannel retailing},
  author={Arslan, Ay{\c{s}}e N and Klibi, Walid and Montreuil, Benoit},
  journal={European Journal of Operational Research},
  volume={294},
  number={3},
  pages={1042--1058},
  year={2021},
  publisher={Elsevier}
}

@article{birge1982value,
  title={The value of the stochastic solution in stochastic linear programs with fixed recourse},
  author={Birge, John R},
  journal={Mathematical programming},
  volume={24},
  number={1},
  pages={314--325},
  year={1982},
  publisher={Springer}
}

@article{escudero2007value,
  title={The value of the stochastic solution in multistage problems},
  author={Escudero, Laureano F and Gar{\'\i}n, Araceli and Merino, Mar{\'\i}a and P{\'e}rez, Gloria},
  journal={Top},
  volume={15},
  number={1},
  pages={48--64},
  year={2007},
  publisher={Springer}
}

@article{wang2024global,
  title={Global multi-sourcing network design with inventory planning under uncertainty},
  author={Wang, Minke and Amiri-Aref, Mehdi and Klibi, Walid and Babai, M Zied},
  journal={International Journal of Production Research},
  pages={1--27},
  year={2024},
  publisher={Taylor \& Francis}
}

@article{yang2023distributionally,
  title={Distributionally robust multi-period location-allocation with multiple resources and capacity levels in humanitarian logistics},
  author={Yang, Yongjian and Yin, Yunqiang and Wang, Dujuan and Ignatius, Joshua and Cheng, TCE and Dhamotharan, Lalitha},
  journal={European Journal of Operational Research},
  volume={305},
  number={3},
  pages={1042--1062},
  year={2023},
  publisher={Elsevier}
}

@techreport{Soundararajan2023,
  author    = {Kirthika Soundararajan},
  year      = {2023},
  title     = {Preventing resource underutilization},
  institution = {Rocketlane},
  type      = {Technical Report}
}

@techreport{LogiNext2025,
  author      = {LogiNext},
  year        = {2025},
  title       = {Pop-up warehouse: Future of delivery management},
  institution = {LogiNext Solutions},
  type        = {Technical Report},
}

@article{correia2018stochastic,
  title={A stochastic multi-period capacitated multiple allocation hub location problem: Formulation and inequalities},
  author={Correia, Isabel and Nickel, Stefan and Saldanha-da-Gama, Francisco},
  journal={Omega},
  volume={74},
  pages={122--134},
  year={2018},
  publisher={Elsevier}
}

@article{kaut2014multi,
  title={Multi-horizon stochastic programming},
  author={Kaut, Michal and Midthun, Kjetil T and Werner, Adrian S and Tomasgard, Asgeir and Hellemo, Lars and Fodstad, Marte},
  journal={Computational Management Science},
  volume={11},
  number={1},
  pages={179--193},
  year={2014},
  publisher={Springer}
}

@article{armour1963heuristic,
  title={A heuristic algorithm and simulation approach to relative location of facilities},
  author={Armour, Gordon C and Buffa, Elwood S},
  journal={Management science},
  volume={9},
  number={2},
  pages={294--309},
  year={1963},
  publisher={Informs}
}

@article{melo2009facility,
  title={Facility location and supply chain management--A review},
  author={Melo, M Teresa and Nickel, Stefan and Saldanha-Da-Gama, Francisco},
  journal={European journal of operational research},
  volume={196},
  number={2},
  pages={401--412},
  year={2009},
  publisher={Elsevier}
}

@article{wesolowsky1973dynamic,
  title={Dynamic facility location},
  author={Wesolowsky, George O},
  journal={Management Science},
  volume={19},
  number={11},
  pages={1241--1248},
  year={1973},
  publisher={INFORMS}
}

@article{klibi2024,
  title={Make smarter investments in resilient supply chains},
  author={Klibi, Walid and Trepte, Kai and Rice Jr, James B},
  journal={MIT Sloan Management Review},
  volume={66},
  number={1},
  pages={53--57},
  year={2024},
  publisher={Massachusetts Institute of Technology, Cambridge, MA}
}

@article{kim2023network,
  title={Network inventory deployment for responsive fulfillment},
  author={Kim, Nayeon and Montreuil, Benoit and Klibi, Walid and Babai, M Zied},
  journal={International Journal of Production Economics},
  volume={255},
  pages={108664},
  year={2023},
  publisher={Elsevier}
}

@incollection{crainic2023hyperconnected,
  title={Hyperconnected city logistics: a conceptual framework},
  author={Crainic, Teodor Gabriel and Klibi, Walid and Montreuil, Benoit},
  booktitle={Handbook on city logistics and urban freight},
  pages={398--421},
  year={2023},
  publisher={Edward Elgar Publishing}
}

@article{jouida2021profit,
  title={Profit maximizing coalitions with shared capacities in distribution networks},
  author={Ben Jouida, Sihem and Guajardo, Mario and Klibi, Walid and Krichen, Saoussen},
  journal={European Journal of Operational Research},
  volume={288},
  number={2},
  pages={480--495},
  year={2021},
  publisher={Elsevier}
}

@article{roels2017win,
  title={Win-win capacity allocation contracts in coproduction and codistribution alliances},
  author={Roels, Guillaume and Tang, Christopher S},
  journal={Management science},
  volume={63},
  number={3},
  pages={861--881},
  year={2017},
  publisher={INFORMS}
}

@article{vanmieghem2003capacity,
  title={Capacity management, investment, and hedging: Review and recent developments},
  author={Van Mieghem, Jan A.},
  journal={Manufacturing \& Service Operations Management},
  volume={5},
  number={4},
  pages={269--302},
  year={2003},
  publisher={INFORMS}
}

@article{jaoui2025multi,
  title={Multi-cycle production network design under supply uncertainty},
  author={Jaoui, Nadia and Klibi, Walid and El Hachemi, Nizar and Aouam, Tarik and Fender, Michel},
  journal={European Journal of Operational Research},
  year={2025},
  publisher={Elsevier}
}

@misc{cbre2025us, 
title = {{U.S. Industrial Figures: Q2 2025}}, 
author = {{CBRE Research}}, 
year = {2025}, 
institution = {CBRE, Inc.}, 
url = {https://www.cbre.com/insights/figures/q2-2025-us-industrial-and-logistics-figures}, 
}

@misc{cbre2025eu,
  title        = {{European Logistics Leasing Figures: Q2 2025}},
  author       = {{CBRE Research}},
  year         = {2025},
  institution  = {CBRE, Inc.},
  url          = {https://www.cbre.com/insights/figures/european-logistics-leasing-figures-q2-2025},

}

@article{singh2009dantzig,
  title={Dantzig-Wolfe decomposition for solving multistage stochastic capacity-planning problems},
  author={Singh, Kavinesh J and Philpott, Andy B and Wood, R Kevin},
  journal={Operations Research},
  volume={57},
  number={5},
  pages={1271--1286},
  year={2009},
  publisher={INFORMS}
}

@article{huang2009value,
  title={The value of multistage stochastic programming in capacity planning under uncertainty},
  author={Huang, Kai and Ahmed, Shabbir},
  journal={Operations Research},
  volume={57},
  number={4},
  pages={893--904},
  year={2009},
  publisher={INFORMS}
}

\pagebreak

\begin{appendices}

\section{Notation of Mathematical Modeling}
\label{appendix_notation}

To formulate the mathematical model in Section \ref{sec:formulation}, we consider the following notation:

{
\begin{longtable}[H]{p{0.06\linewidth} p{0.92\linewidth}} 
\caption{Notation}\\
\hline
Sets \\
\hline
$\boldsymbol{I}$ 
& set of demand locations indexed by $i$\\
$\boldsymbol{J}$ 
& set of potential facility locations indexed by $j, j^{'}$\\
$\boldsymbol{IJ}$ 
& set of feasible demand-facility pairs (under distance limit threshold)\\
$\boldsymbol{JJ}$ 
& set of feasible facility-facility capacity module relocation pairs (new rental $0J$ \& return $J0$)\\
$\boldsymbol{L}$ 
& set of capacity levels indexed by $l$\\
\hline
\multicolumn{2}{l}{Parameters}\\
\hline
$w_{t}$ 
& binary parameter, whether to change the capacity level in stage $t$\\
$u_{jl}$ 
& number of capacity modules in facility $j$ at capacity level $l$\\
$v_j$ 
& initial capacity level at facility $j$\\
$p_n$ 
& probability of node $n$\\
$k_{jn}$ 
& average throughput rate per capacity module at facility $j$ at node $n$\\
$d_{in}$ 
& demand volume in demand location $i$ at node $n$\\
$c_{j l_1 l_2 n}^F$ 
& combined cost to change the capacity level of facility $j$ from $l_1$ to $l_2$ at node $n$\\
$c_{j j^{'} n}^A$ 
& combined allocation and relocation cost per capacity module from facility $j$ to $j^{'}$ at node $n$\\
$c_{i j n}^T$ 
& transportation cost per unit of demand from facility $j$ to demand location $i$ at node $n$\\
$c_{n}^O$ 
& outsourcing cost per unit of demand at node $n$\\
\hline
\multicolumn{2}{l}{Decision variables}\\
\hline
$Y_{j l_1 l_2 n}$ 
& binary variable, whether to change the capacity level of facility $j$ from $l_1$ to $l_2$ at node $n$\\
$F_{jj^{'}n}$ 
& integer variable, number of capacity modules transported from facility $j$ to $j^{'}$ at node $n$\\
$S_{jn}$ 
& integer variable, number of capacity modules available in facility $j$ at node $n$\\
$X_{ijn}$ 
& continuous variable, demand volume in demand location $i$ assigned to facility $j$ at node $n$\\
$R_{jn}$ 
& continuous variable, outsourced demand volume in facility $j$ at node $n$\\
\hline
\vspace{-0.05in}
\label{tab:notation}
\end{longtable}
}

We note that the combined cost $c^F_{j l_1 l_2}$ represents multiple costs of facility opening, closing, reopening, and capacity adjustment, while $c^A_{j j^{'}}$ denotes the costs to expand, reduce, and relocate capacity modules within the facility network. More details about these cost parameters can be found in Appendix \ref{appendix_b}.

\section{Enhanced SDDiP with Alternating Cut Strategy}
\label{appendix_algo_1}

The alternating cut strategy was first introduced by \cite{angulo2016improving} for the two-stage integer L-shaped method, where the algorithm alternates between solving linear and mixed-integer subproblems to improve computational performance. More recently, \cite{arslan2021distribution} successfully applied the integer L-shaped method to address a complex two-stage stochastic model for omnichannel distribution network design. To the best of our knowledge, our work is the first to extend and integrate this alternating cut concept within the multi-stage SDDiP framework. Specifically, in the backward pass of SDDiP, we tailor an alternating cut strategy to balance the contribution of LP-based and integer-based cuts. In our implementation, LP-based cuts are generated by solving the relaxation (e.g., LP relaxation and Lagrangian relaxation) of the forward subproblems, which allows them to be computed quickly but without guaranteeing tightness. Whereas integer-based cuts are generated by solving the mixed integer programs, and thus usually provide stronger bounds but at a higher computational cost. Among the cut types listed in Sections \ref{sec:cg} and \ref{sec:scg}, all cuts except the integer optimality cut and the Lagrangian cut fall into the LP-based category. 

The proposed SDDiP enhancement with alternating cut strategy is summarized in Algorithm \ref{alg:opt_alg}. Specifically, for each state $(t,n,y^*)$, we maintain two memory sets, denoted as $V_{LP}$ and $V$. When a state is encountered for the first time, the algorithm generates an LP-based cut from the relaxed subproblem and stores this state in $V_{LP}$. If the same state is revisited in later iterations, the algorithm instead generates an integer-based cut and records it in $V$. To the best of our knowledge, this is the first work to adapt the alternating cut strategy within the SDDiP framework. By systematically alternating between LP-based and integer-based cuts, our enhanced method accelerates convergence by leveraging global relaxations and tighter integer-based approximations in the backward pass. 

\begin{algorithm}[!htbp]
{
\renewcommand{\baselinestretch}{1.1}\selectfont
\caption{{Enhanced SDDiP with Alternating Cut Strategy}}
\KwData{Initial lower-approximation cut set $\Phi_t, \forall t \in \boldsymbol{T}$, strategy threshold $\zeta$} 
Iteration number $k \gets 1$, accept number $a \gets 0$;  Memory sets $V_{LP} \gets \varnothing$, $V \gets \varnothing$\;
Lower bound $LB^k \gets -\infty$, Upper bound $UB^k \gets +\infty$\;
\While{stopping criterion is not satisfied}{
  
  Randomly sample $M$ scenarios $\boldsymbol{\Omega^k} = \{ \boldsymbol{\omega^k_{1}}, ..., \boldsymbol{\omega^k_{m}}, ..., \boldsymbol{\omega^k_{M}}\}$\;
  \For{$m = 1,...,M$}{ 
    \For{$t = 1,...,T$}{
        Solve forward subproblem $\underline{Q}^k_t( \hat{X}^{k}_{t-1}, \Phi_t^k, \boldsymbol{\omega^{k}_{t m}})$\; 
        Record the optimal solution $\hat{X}^k_{t m}$, optimal objective value $z^k_{t m}$, and the value of lower-approximation $\hat{\theta}_{tm}^k$\;    
    }
    Compute modified objective value $z^k_{m} \gets \sum_{t \in \boldsymbol{T}} \left( z^k_{t m} + \hat{\theta}^k_{tm} \right)$\;
  }
  Calculate mean $\bar{z}^k \gets \frac{1}{M} \sum_{m=1}^{M} z^k_m$ and variance $(\sigma^k)^2 \gets \frac{1}{M-1} \sum_{m=1}^{M} (z^k_m - \bar{z}^k)^2$\;
  Update statistical upper bound $UB^k \gets \bar{z}^k + \mathit{z}_{\alpha/2} \cdot \frac{\sigma^k}{\sqrt{M}}$\;


  \uIf{$a >= \zeta$}{
  Memory sets $V_{LP} \gets \varnothing$, $V \gets \varnothing$; Accept number $a \gets 0$\;
  }
  \For{$t=T,\dots,2$}{
    \For{sampled node $n$ at stage $t$}{
      \For{child node $c$ of $n$}{    
        Extract state decisions $y^*=(\hat{y}^t_n)$ from forward solution $\hat{X}^{k}_{t}$\;
        Create a unique identifier $state\_id \leftarrow (t,n,y^*)$\;  
        \uIf{$state\_id \notin V_{LP}$}{
          Solve LP relaxation backward subproblem at child node $c$ with fixed $y^*$\;
          Generate and add LP-based cut to $\Phi_t^{k+1}$\;
        $V_{LP} \leftarrow V_{LP} \cup \{state\_id\}$\;
        }
        \uElseIf{$state\_id \notin V$}{
          Solve integer backward subproblem at child node $c$ with fixed $y^*$\;
          Generate an integer optimality cut or Lagrangian cut to $\Phi_t^{k+1}$\;
          $V \leftarrow V \cup \{state\_id\}$\;
        }
        \Else{
          Accept the node and no cut added; $a \gets a + 1$\;
        }
      }
    }
  }
  Solve $\underline{Q}_{1}^{k}(v, \Phi^{k+1}_{1})$, where $v$ is the initial facility capacity information\;
  Record optimal solution $\hat{X}^k_{1}$ and optimal value $z^k_1$; Update lower bound $LB^k \gets z^k_1$\;
  $k \gets k+1$ \;
} 
\label{alg:opt_alg}
}
\end{algorithm}

In the original alternating cut strategy for Benders decomposition, once a node is accepted, it is permanently marked and no further cuts are generated from it. However, in the context of SDDiP, the situation is different: even if a state $(t,n,y^*)$ has been used to generate both LP-based and integer-based cuts in one iteration, revisiting it in a later iteration may still yield new, non-redundant cuts. This is because the dual information and incumbent solutions evolve as the forward and backward passs progress, meaning that previously visited states can provide new cut coefficients reflecting updated value function approximations. To accommodate this, we periodically clear the sets 
$V_{LP}$ and $V$, thereby allowing states to re-enter the alternating process in subsequent iterations. This reset mechanism prevents states from being prematurely excluded and ensures that the algorithm continues to exploit new structural information over time.

The frequency of memory-set clearing is controlled by a threshold parameter $\zeta$, defined in terms of the number of accepted states. When the number of accepts $a$ reaches $\zeta$, both $V_{LP}$ and $V$ are cleared, and the acceptance counter $a$ is reset to zero. This has a direct effect on the balance between LP-based and integer-based cuts. If the threshold is set small, the memory sets are cleared more often, meaning that many states will be encountered as “first-time” visits, leading to more LP-based cuts. Conversely, when the threshold is larger, states tend to persist longer in the memory, and the algorithm generates more integer-based cuts when revisited. In this way, the parameter $\zeta$ naturally regulates the mix of cutting planes: smaller values promote broader exploration of relaxations of forward subproblems, while larger values bias the algorithm toward stronger but more computationally expensive integer cuts. This balance helps guide convergence by preventing dominance of either cut type and by allowing the convergence of the algorithm to benefit from mixed cutting planes.

\section{SDDiP Parallelization Framework}
\label{appendix_a}

In the SDDiP algorithm, the need to solve numerous subproblems at each stage in both the forward and backward passs presents a potential opportunity for acceleration through parallelization. While the concept is not new, existing literature typically focuses on limited implementations. For example, \cite{zou2018multistage} applied parallelization only to the backward pass of SDDiP, accelerating cut generation but leaving the forward pass sequential. \cite{lara2020electric} extended this by proposing a parallelization of both forward and backward passs for electric power infrastructure planning, but their implementation is limited to a three-processor structure. To fill this gap, our work builds upon these by developing a more general and flexible parallelization framework for SDDiP. The proposed framework is designed to be scalable for an arbitrary number of processes and parallelizes the subproblem-solving tasks in both the forward and backward passs, enabling efficient cut generation and aggregation. 

Algorithm \ref{alg:parallel_sddip} describes the proposed parallelization framework. In each iteration, the procedure begins with a preprocessing phase, where one master process ($p=1$) samples a subset of nodes and distributes them across available processors. Each processor then executes the forward and backward passs in parallel. In the forward pass, stage-node subproblems are solved independently according to the preprocessing assignment. After each stage, synchronization barriers are placed to ensure that all subproblems of a given stage are completed before any process moves to the next. The master process then updates the upper bound. In the backward pass, each processor solves its assigned subproblems to compute cut coefficients. After another stage-wise synchronization, these coefficients are aggregated to form and add new cuts to the lower-approximation cut sets $\Phi_t$. The master process subsequently solves the origin node at stage $1$ to update the lower bound.

\begin{algorithm}[!htbp]
{
\renewcommand{\baselinestretch}{1.1}\selectfont
\caption{{SDDiP Parallelization Framework}}
\KwData{Initial lower-approximation cut set $\Phi_t$, number of processes $P$, number of sampled nodes per iteration $m$} 
\While{stopping criterion is not satisfied}{
  
    \If{$p=1$}{
      Randomly sampling $m$ nodes\;
      Distribute sampled nodes to each process $p=1,\dots,P$\;
      Generate process-node assignment set $forward\_nodes$ and $backward\_nodes$\;
    }
    barrier()\;
    
    \ForEach{process $p=1,\dots,P$ \textbf{in parallel}}{
        
        \If{$p=1$}{
          Solve forward subproblem for origin node at stage $1$\;
        }
        barrier()\;
        \For{$t=2,\dots,T$}{
          \For{node $n$ in forward\_nodes[$t,p$]}{
            Solve forward subproblem for node $n$ at stage $t$\;
          }
          barrier()\;
        }
        \If{$p=1$}{
          Update upper bound $UB^k$\;
        }
        barrier()\;

        
        \For{$t=T,\dots,2$}{
          \For{node $n$ in backward\_nodes[$t,p$]}{
            Solve backward subproblem for node $n$ at stage $t$ to collect cut coefficients\;
          }
          barrier()\;
          Add cuts in $\Phi_t$ using collected coefficients\;
          barrier()\;
        }
        \If{$p=1$}{
          Solve backward subproblem for origin node at stage $1$\;
          Update lower bound $LB^k$\;
        }
        barrier()\;
    }
  }
\label{alg:parallel_sddip}
}
\end{algorithm}

Specifically, to enable balanced workload distribution in the preprocessing step, we partition the set of nodes as evenly as possible among the processes. At the beginning of each iteration, the master process samples nodes from the scenario tree for every stage. These sampled nodes form the subproblems that need to be solved at that stage. For the forward pass, the original node at stage $1$ is always assigned to the master process, while nodes at later stages are evenly partitioned among processes. Partitioning is done by dividing the set of nodes into nearly equal subsets so that each process receives a comparable number of nodes, with remainders distributed in a balanced way. For each stage in the backward pass, child nodes are first generated from the sampled parent nodes of that stage, and the resulting set of child nodes is partitioned among processes using equal distribution. This approach ensures that the workload at each stage of the scenario tree is allocated fairly across processors, preventing inefficiency caused by uneven node assignment. Moreover, we remark that the implementation of this framework is developed in Python using the \texttt{multiprocessing} library. To manage inter-process communication and data consistency, we utilize shared dictionaries \texttt{manager.dict()}. These dictionaries record critical information, including state-variable forward solutions, backward cut coefficients, and process-node assignments, which can be accessed and updated across all processes.

\section{Complementary Case Study Parameters and Setup}
\label{appendix_b}

\subsection{Modular Construction Context}
\label{app:modular}

Modular construction offers a faster, cheaper, and greener alternative to traditional site-built methods. Building sections or modules are prefabricated in factory settings and transported to erection sites for installation, enabling parallel processing of site preparation and building fabrication. This leads to faster completion times, reduced on-site labor, and significant cost savings compared with site-built construction.

A key challenge in modular construction is the transportation of modules. Road restrictions, expensive transportation costs, and damage risk during transit associated with shipping oversize modules make finished housing modules highly sensitive to shipping distances, reinforcing the importance of an efficient deployment of assembly centers.

\subsection{Region Instances and Geographic Scope}
\label{app:regions}

To evaluate scalability, we define three nested regional instances based on the partner company's project and site data.

\subsubsection*{SE Instance (Southeast).}
The SE instance represents a small region covering the southeastern United States and includes:
\begin{itemize}
    \item 7 potential facility locations
    \item 50 project sites
    \item Total demand: 5{,}172 capacity modules across the planning horizon
\end{itemize}

\subsubsection*{EAST Instance.}
The EAST instance roughly doubles the scale of the SE instance and includes:
\begin{itemize}
    \item 14 potential facility locations
    \item 81 project sites
    \item Total demand: 8{,}020 capacity modules
\end{itemize}

\subsubsection*{US Instance.}
The US instance represents a large-scale nationwide network and includes:
\begin{itemize}
    \item 28 potential facility locations
    \item 131 project sites
    \item Total demand: 13{,}520 capacity modules
\end{itemize}

Candidate locations (red points in Figure~\ref{fig:network_info} (b)) reflect real sites identified by the company as feasible for leasing assembly centers, while project locations (blue points in Figure~\ref{fig:network_info} (a)) represent upcoming projects through 2026, with point size indicating forecasted module requirements.

\subsection{Planning Horizon and Time Structure}
\label{app:horizon}

The planning horizon is set to a maximum of 12 months, with each operational period corresponding to one month.To assess scalability, we also consider shorter horizons of 3, 6, and 9 months, which directly change the depth and size of the scenario tree and thus overall problem complexity.

\subsection{Capacity Levels}
\label{app:capacity}

We test two capacity structures per facility to represent different degrees of modularity.

\subsubsection*{2-Capacity-Level Variant.}

In the 2-level variant, each facility can operate at two distinct capacity levels (e.g., Level 1: 0 modules, Level 2: 6 modules), capturing facilities with limited modular flexibility.

\subsubsection*{3-Capacity-Level Variant.}

In the 3-level variant, each facility can operate at three distinct capacity levels (e.g., Level 1: 0 modules, Level 2: 3 modules, Level 3: 6 modules), providing greater flexibility to match capacity to realized demand.

\subsection{Cost Parameters}
\label{app:costs}

Cost parameters are calibrated from the partner company's internal accounting data and reflect typical rates in modular construction projects.

\subsubsection*{Facility Capacity Adjustment Costs.}

These combined costs include:
\begin{itemize}
    \item Facility opening cost (from closed to active)
    \item Facility closing cost (from active to closed)
    \item Facility reopening cost (from closed back to active)
    \item Capacity level transition costs for upgrading or downgrading the number of deployed modules
\end{itemize}

\subsubsection*{Module Relocation and Deployment Costs.}

Relocation and deployment costs cover:
\begin{itemize}
    \item Transportation cost to relocate modules between facilities or between depot and facilities
    \item Handling and setup costs at receiving facilities
    \item Costs associated with returning modules to the depot
\end{itemize}

\subsubsection*{Demand Fulfillment Costs.}

Demand fulfillment costs are modeled as:
\begin{itemize}
    \item Transportation cost per unit of demand from each facility to each demand location
    \item Outsourcing cost per unit of demand, representing payments to external providers or expensive on-demand capacity
\end{itemize}

\subsection{Uncertainty Modeling}
\label{app:uncertainty}

\subsubsection*{Demand Uncertainty.}

The company provides monthly demand forecasts for each project location through 2026; actual monthly demands are modeled as stochastic realizations around these forecasts. Demand scenarios are generated by:
\begin{itemize}
    \item Drawing realizations from probability distributions calibrated on historical project data
    \item Assuming stage-wise independence conditional on the forecast
    \item Capturing both upside and downside deviations relative to baseline forecasts
\end{itemize}

\subsubsection*{Supply Disruption Uncertainty.}

Module throughput at each facility may be reduced by disruptions such as equipment failures or labor shortages. Supply scenarios are modeled by:
\begin{itemize}
    \item Assigning disruption probabilities to facility–month pairs
    \item Reducing effective throughput to a fraction of nominal capacity (e.g., 70–80\%) in disruption states
    \item Combining demand and supply realizations into a joint scenario tree
\end{itemize}

\subsection{Scenario Tree Construction}
\label{app:scenario-tree}

The scenario tree captures the joint evolution of demand and supply uncertainties over the planning horizon.

\subsubsection*{Branching Structure.}

At each stage (month), the tree branches to represent possible combinations of demand and supply realizations, with branching factors 4 per stage. For a branching factor \(b\) and horizon \(T\), the number of leaf nodes scales as \(b^T\), motivating additional sampling within SDDiP.

\subsubsection*{Monte Carlo Scenario Sampling.}

To maintain tractability, we embed Monte Carlo sampling in SDDiP:
\begin{itemize}
    \item In the forward pass, we sample \(M = 5\) scenarios from the tree
    \item In the backward pass, cuts are generated using the full tree
    \item This approach reduces per-iteration cost while progressively exploring the scenario space
\end{itemize}

\subsubsection*{Scenario Generation Data.}

Scenario generation relies on:
\begin{itemize}
    \item Three years of historical monthly demand data from the partner company
    \item Empirical disruption rates observed across facilities
    \item Probability distributions fitted to demand variability and disruption frequency
\end{itemize}

\subsection{Implementation Details}
\label{app:implementation}

\subsubsection*{Optimization Solver.}

All subproblems in SDDiP are solved with Gurobi 11.0.3 as a mixed-integer programming solver.

\subsubsection*{Computing Environment.}

Experiments are run on a 24-core Intel Xeon Gold 6226 machine with 128 GB RAM using Python 3.9.12.

\subsubsection*{Parallelization.}

The parallelization framework described in Section~4.3.3 is implemented using Python's \texttt{multiprocessing} library with shared dictionaries to synchronize data across processes. The number of processes is set to 4 in the experiments of this study.

\subsubsection*{Convergence Criteria.}

The SDDiP algorithm terminates when either:
\begin{itemize}
    \item The optimality gap \(\frac{UB^k - LB^k}{UB^k} \leq 1\%\) is reached, or
    \item A predefined runtime limit (e.g., 24 hours for large instances) is exceeded
\end{itemize}

\subsubsection*{Sampling and Cut Management.}

Sampling and cut management are configured as follows:
\begin{itemize}
    \item Forward pass: \(M = 5\) scenarios per iteration
    \item Backward pass: all nodes in the scenario tree are visited to generate cuts
    \item Cut aggregation: cuts from nodes at the same stage are aggregated to update the lower-approximation functions \(\phi_t(Y_t)\)
\end{itemize}

\subsection{Cost Parameter Estimation}
\label{appendix_b_1}

\vspace{0.1in}
{
\renewcommand\labelitemi{\tiny$\bullet$}
\noindent
\textbf{Parameters}
\begin{itemize}
    \setlength\itemsep{0.12em}
    \item $c_{j l_1 l_2 n}^F$: combined cost to change the capacity level of facility $j$ from $l_1$ to $l_2$
    \item $c_{jl}^C$: commissioning cost for opening or reopening facility $j$ at capacity level $l$
    \item $c_{jl}^D$: decommissioning cost for closing facility $j$ at capacity level $l$
    \item $c_{jl}^R$: cost to reduce the capacity level of active facility $j$ by $l$
    \item $c_{jl}^E$: cost to expand the capacity level of active facility $j$ by $l$
    \item $c_{jl}^M$: cost to maintain the capacity level of active facility $j$ at $l$
    \item $c_{j j^{'} n}^A$: combined allocation and relocation cost per production unit from facility $j$ to facility $j^{'}$
    \item $c^P$: commissioning cost to rent a new production unit
    \item $c^S$: decommissioning cost to return an unused production unit
    \item $c_{jj^{'}}^U$: cost to transport a production unit from facility $j$ to facility $j^{'}$
\end{itemize}
}

\vspace{0.1in}
We note that the combined cost $c^F_{j l_1 l_2}$ can represent multiple costs of facility opening, closing, reopening, and capacity changes, while $c^A_{j j^{'}}$ denotes costs to purchase, store, and relocate capacity modules: 

\vspace{0.1in}
\begin{equation}
\begin{aligned}
c^F_{j l_1 l_2}
&=\left\{\begin{array}{cl}
c^C_{j l_2} + c^M_{j l_2}, & \text { if } l_2 > l_1 = 0 \text { (open facility) }\\
c^D_{j l_1}, & \text { if } l_1 > l_2 = 0 \text { (close facility) }\\
c^E_{j, l_2-l_1} + c^M_{j l_2}, & \text { if } l_2 > l_1 > 0 \text { (expand cap) }\\
c^R_{j, l_1-l_2} + c^M_{j l_2}, & \text { if } l_1 > l_2 > 0 \text { (reduce cap) }\\
c^M_{j l_2}, & \text { if } l_1 = l_2 \text { (maintain cap) }\\
\end{array} \right.
\\
&\qquad \qquad \qquad \qquad \qquad \qquad  j \in \boldsymbol{J}, \; l_1, l_2 \in \boldsymbol{L}.
\end{aligned}
\end{equation}

\vspace{0.1in}
\begin{equation}
\begin{aligned}
c^A_{j j^{'}}&=\left\{\begin{array}{cl}
c^P, & \text { if } j=0 \text { (new rent) }\\
c^S, & \text { if } j^{'} = 0 \text { (return) }\\
c_{j j^{'}}^U, & \text { if } j \neq 0 \text { and } j^{'} \neq 0 \text { (relocation) }\\
\end{array} \right.
\\
&\qquad \qquad \qquad \qquad \qquad \qquad j,j^{'} \in \boldsymbol{J} \cup \{0\}.
\end{aligned}
\end{equation}

\vspace{0.05in}
\subsection{Input Parameter Summary}
\label{appendix_b_2}
Table \ref{tab:other_paras} summarizes the key parameters used in this study.
{
\begin{table}[h!]
    \begin{center}
    \caption{{Input parameter summary.}}
    \resizebox{0.76\linewidth}{!}{
    \begin{tabular}{cll}
    Notation & Implication & Value\\
    \midrule
    $C^{C}$ & unit commissioning cost for opening or reopening facility & \$50,000 \\
    $C^{D}$ & unit decommissioning cost for closing facility & \$25,000 \\
    $C^{E}$ & unit cost to expand one capacity level & \$12,500 \\
    $C^{R}$ & unit cost to reduce one capacity level & \$6,250 \\
    $C^{M}$ & unit cost to maintain one capacity level & \$22,188 \\
    $C^P$ & commissioning cost to rent a new production unit & \$19,104\\
    $C^S$ & decommissioning cost to return an unused production unit & \$19,104 \\
    $C^U$ & transportation cost per production unit per hour & \$60/hour \\
    $C^T$ & transportation cost per housing module per hour & \$120/hour \\
    $C^O$ & outsourcing penalty per housing module & \$10,000\\
    $L^M$ & distance limit of housing module transportation & 250 miles\\
    $L^P$ & distance limit of production unit transportation & 750 miles\\
    \midrule
    \label{tab:other_paras}
    \end{tabular}
    }
    \vspace{-0.3in}
    \end{center}
\end{table}
}

\section{Scenario Tree Generation}
\label{appendix_scenario_tree_generation}

With the assumption of stage-wise independence, the evolution of the system at each stage is independent of past realizations, so that uncertainty at stage $t$ can be modeled without conditioning on previous stages beyond the current state variables. To construct such a recombining scenario tree, we begin with a given demand $D_{i1}$ in the first stage for project $i \in \boldsymbol{I}$. In following stages $t=2,\dots,T$, stochastic demand realizations $\tilde{D}_{it}$ are generated from a nominal demand $D_{it}$ as below: 
\begin{equation}
\tilde{D}_{it} = D_{it} \cdot \xi_t, \quad \xi_t \sim TN(1, (t-1)^2\sigma^2),
\end{equation}

where $\xi_t$ is a truncated normal random variable with mean $1$ and variance increasing quadratically in $t-1$.The nominal demand $\tilde{D}_{it}$ is given by the number of housing modules predicted by our partner company. The demand standard deviation parameter $\sigma$ controls the intensity of demand uncertainty, where a higher $\sigma$ corresponds to larger deviations from the nominal demand forecast. This construction captures the pattern that forecast errors typically accumulate and magnify over longer planning horizons. 

For production rate uncertainty, we model supply-side disruptions using a Poisson distribution. Specifically, the number of disruptions in facility $j$ during stage $t$ is drawn from $x_{jt} \sim \mathcal{P}(\lambda)$, where the disruption rate $\lambda$ controls the intensity of supply uncertainty. Then, the effective production rate is generated from the nominal production rate $P_{jt}$ according to:
\begin{equation}
\tilde{P}_{jt} = P_{jt} \cdot \delta_{jt}, \quad \delta_{jt} \sim U(a, b),
\end{equation}

In our setting, the nominal rate of a production unit is two modules per day based on realistic factors. Accordingly, $P_{jt}$ at facility $j$ in period $t$ is calculated as the product of the number of production units assigned to that facility and the nominal rate of one unit. To capture uncertainty, the stochastic factor $\delta_{jt}$ is modeled as a uniform random variable on $[a, b]$ with $a, b \in [0,1]$, representing the proportion of the nominal rate that is actually realized. The parameters $a$ and $b$ represent the intensity of production rate reduction and are determined by the number of disruptions $x_{jt}$ through the following relationship:

{
\begin{equation}
(a, b) =
\begin{cases}
(1.0, 1.0) & \text{if } x_{jt} = 0 \text{  (no disrupt)} \\
(0.8, 0.99) & \text{if } x_{jt} = 1 \text{  (small disrupt)} \\
(0.6, 0.79) & \text{if } x_{jt} = 2 \text{  (med disrupt)} \\
(0, 0.59) & \text{if } x_{jt} \geq 3 \text{  (large disrupt)}.
\end{cases}
\end{equation}    
}

In this formula, as the number of disruptions $x_{jt}$ increases, the potential reduction in the production rate $\delta_{jt}$ becomes more severe, indicating a greater operational impact. We note that, in the following experiments, both the demand standard deviation $\sigma$ and disruption rate $\lambda$ are fixed at their default value of 0.5. While in Section \ref{sec:vss}, we explore model performance under different uncertainty levels by altering these two parameters. Each parameter setting leads to a distinct scenario tree, and during the SDDiP forward pass, we employ classical Monte Carlo sampling with four sampled scenarios per iteration to approximate the stochastic process.

\section{Problem Size of Test Instances}
\label{appendix_size}

Table \ref{tab:modelsize} presents the detailed problem sizes of all test instances used in this study.

{
\begin{table}[h!]
\begin{center}
\caption{{Model size by capacity levels, region, and planning horizon.}}
\label{tab:modelsize}
\resizebox{0.93\textwidth}{!}{%
\begin{tabular}{cccccccc}
\toprule
Capacity levels & Region & \#Periods & \#Continuous variables & \#Integer variables & \#Constraints & \#Non-zeros \\
\midrule
2 & SE   & 3  & 1,722    & 1,788    & 2,229     & 8,445 \\
 &      & 6  & 111,930  & 116,031  & 144,696   & 551,424 \\
 &      & 9  & 7,165,242 & 7,427,394 & 9,262,395 & 35,301,891 \\
 & EAST & 3  & 2,856    & 3,993    & 4,056     & 16,551 \\
 &      & 6  & 185,640  & 259,356  & 263,451   & 1,081,002 \\
 &      & 9  & 11,883,816 & 16,602,399 & 16,864,542 & 69,205,677 \\
 & US   & 3  & 4,347    & 7,311    & 7,458     & 29,025 \\
 &      & 6  & 282,555  & 475,026  & 484,581   & 1,897,188 \\
 &      & 9  & 18,087,867 & 30,408,597 & 31,020,264 & 121,459,431 \\
\midrule
3 & SE   & 3  & 1,722    & 2,523     & 2,964     & 12,939 \\
 &      & 6  & 111,930  & 163,806   & 192,471   & 847,566 \\
 &      & 9  & 7,165,242 & 10,485,729 & 12,320,730 & 54,263,505 \\
 & EAST & 3  & 2,856    & 5,463     & 5,526     & 25,539 \\
 &      & 6  & 185,640  & 354,906   & 359,001   & 1,673,286 \\
 &      & 9  & 11,883,816 & 22,719,069 & 22,981,212 & 107,128,905 \\
 & US   & 3  & 4,347    & 10,251    & 10,398    & 47,001 \\
 &      & 6  & 282,555  & 666,126   & 675,681   & 3,081,756 \\
 &      & 9  & 18,087,867 & 42,641,937 & 43,253,604 & 197,305,887 \\
\bottomrule
\end{tabular}
}
\end{center}
\end{table}
}

\section{Complementary Experimental Results}

\subsection{Computational Performance Comparison across Different Cut Combinations}
\label{appendix_numerical_1}

For the first set of experiments in Section \ref{sec:nr}, we focus on identifying the most effective combination of cutting planes used in the SDDiP backward pass. Among these cut types, the integer optimality cuts (I), Lagrangian cuts (L), Benders cuts (B), and strengthened Benders cuts (SB) were initially applied to SDDiP in \cite{zou2019stochastic}, which serve as our benchmark. We compare these against the new set of cuts introduced in this paper: the Pareto-optimal cut (PT), the independent Magnanti-Wong cut (IM), the strengthened Pareto-optimal cut (SPT), and the strengthened independent Magnanti-Wong cut (SIM). Following a similar experimental setting in \cite{zou2019stochastic}, we ensure that each combination includes at least one type of tight cuts (i.e., I and L). The algorithm terminates once the lower bounds obtained in the backward passs stabilize or the current optimality gap falls below 1\%. In addition, in this experiment, we set a time limit of 2 hours. After termination, we independently simulate 500 forward paths to construct a 95\% confidence interval and report its right endpoint as the statistical upper bound of the optimal value. Since the upper bound computed during the forward pass relies on in-sample scenarios, it may be upward-biased. The independent simulation instead provides an out-of-sample estimate, offering a more accurate and statistically reliable measure of convergence in the SDDiP algorithm. To validate the solution quality of SDDiP, we further solve each instance directly using Gurobi to obtain the optimal value and record the true optimality gap in Table \ref{tab:exp1.1} and Table \ref{tab:exp1_2cap}. The performance of each cut combination is evaluated based on five key metrics: the total number of iterations to algorithm termination (\#Iter), the true optimality gap percentage (\%Gap), the total number of backward cuts generated (\#Cuts), the total runtime in seconds (Runtime), and the average runtime per iteration in seconds (Runtime/iter). Finally, we note that this first experiment is performed on small- to medium-scale, focusing on the SE region with a 6-month planning horizon and examining both two and three capacity levels to assess the performance of each cut type. After identifying the most effective cutting planes, we then test the algorithm on large instances in the second set of experiments.

{
\renewcommand{\baselinestretch}{1.3}\selectfont
\begin{table}[h!]
\begin{center}
\caption{{Computational performance comparison across different cut combinations (2 capacity levels).}}
\label{tab:exp1_2cap}
\resizebox{\linewidth}{!}{
\begin{tabular}{cccccccccccc}
\toprule
 \multicolumn{6}{c}{Classical SDDiP cutting planes} & \multicolumn{6}{c}{Proposed SDDiP cutting planes (Section \ref{sec:scg})} \\
\cmidrule(lr){1-6} \cmidrule(lr){7-12}
Cut type & \#Iter & \%Gap & \#Cuts & Runtime & Runtime/iter
& Cut type & \#Iter & \%Gap & \#Cuts & Runtime & Runtime/iter \\
\midrule
I      & 481 & 40.99 & 8,759 & 7,932 & 16.49
& PT + I   & 72 & 0.01 & 2,638 & 7,470 & 103.76 \\
L      & 33  & 0.01  &   613 &   530 & 16.06
& PT + L   & 25 & 0.01 &   932 & 1,649 & 65.94 \\
B + I  & 141 & 0.01  & 5,154 & 2,878 & 20.41
& IM + I   & 81 & 0.01 & 2,960 & 1,183 & 14.61 \\
B + L  & 27  & 0.01  & 1,006 &   575 & 21.29
& IM + L   & 26 & 0.01 &   970 &   569 & 21.88 \\
SB + I & 137 & 0.01  & 5,006 & 3,965 & 28.94
& SPT + I  & 18 & 0.01 &   664 &   679 & 37.72 \\
SB + L & 24  & 0.01  &   894 &   562 & 23.41
& SPT + L  & 15 & 0.01 &   552 &   733 & 48.85 \\
        &     &       &       &       &
& SIM + I  & 28 & 0.01 & 1,044 &   381 & 13.60 \\
        &     &       &       &       &
& SIM + L  & 17 & 0.01 &   626 &   387 & 22.77 \\
\bottomrule
\end{tabular}
}
\end{center}
\vspace*{-6pt}
\end{table}
}

Table~\ref{tab:exp1_2cap} highlights the significant performance improvement of the newly proposed cutting planes compared to classical cutting planes for the 2-capacity-level instances. Specifically, the benchmark I cuts alone fail to achieve convergence, with optimality gaps exceeding 40\%; While, for this small instance, all the other cutting planes achieve the near-optimal solution within time limit. The introduction of B cuts and SB cuts accelerates the convergence speed to close the optimality gap. However, the proposed cuts demonstrate superior performance: the strengthened variants SPT and SIM solve the instance with fewer iterations. Overall, the combination of SIM and I cuts proves to be the most efficient strategy, attaining a 0.01\% optimality gap with the shortest runtime.

\subsection{Computational Performance Comparison across Four SDDiP Variants }
\label{appendix_numerical_2} 

{
\renewcommand{\baselinestretch}{1.25}\selectfont
\begin{table}[h!]
\begin{center}
\caption{{Computational performance of SDDiP algorithm variants for two-capacity-level instances.}}
\label{tab:exp1.2.1}
\resizebox{\linewidth}{!}{
\begin{tabular}{ccccccccccc}
\toprule
 & & & \multicolumn{4}{c}{Classic SDDiP} & \multicolumn{4}{c}{Parallel SDDiP} \\
\cmidrule(lr){4-7} \cmidrule(lr){8-11}
Region & \#Periods & \#Sce & \#Iter & \%Gap & Runtime & Runtime/Iter & \#Iter & \%Gap & Runtime & Runtime/Iter \\
\midrule
SE & 3 & 16 & 14 & 0.02 & 43 & 3.09 & 13 & 0.02 & 33 & 2.56 \\
& 6 & 1,024 & 28 & 0.01 & 409 & 14.60 & 19 & 0.01 & 205 & 10.80 \\
& 9 & 65,536 & 27 & 0.13 & 766 & 28.39 & 30 & 0.13 & 746 & 24.87 \\
& 12 & 4,194,304 & 51 & 0.21 & 3,108 & 60.93 & 39 & 0.21 & 2,147 & 55.06 \\
EAST & 3 & 16 & 25 & 0.45 & 97 & 3.89 & 26 & 0.45 & 79 & 3.03 \\
& 6 & 1,024 & 64 & 0.12 & 1,909 & 29.83 & 74 & 0.12 & 1,510 & 20.41 \\
& 9 & 65,536 & 121 & 0.07 & 12,008 & 99.24 & 195 & 0.07 & 20,123 & 103.20 \\
& 12 & 4,194,304 & 205 & 0.35 & 52,951 & 258.30 & 168 & 0.35 & 30,300 & 180.36 \\
US & 3 & 16 & 42 & 0.05 & 259 & 6.16 & 60 & 0.05 & 292 & 4.87 \\
& 6 & 1,024 & 400 & 5.74 & 85,502 & 213.76 & 342 & 0.11 & 34,333 & 100.39 \\
& 9 & 65,536 & 251 & 13.94 & 91,487 & 364.49 & 308 & 5.33 & 92,487 & 300.28 \\
& 12 & 4,194,304 & 199 & 36.80 & 92,280 & 463.72 & 234 & 36.89 & 93,910 & 401.32 \\
\midrule
 & & & \multicolumn{4}{c}{SDDiP + alternating cut strategy} & \multicolumn{4}{c}{Parallel SDDiP + alternating cut strategy} \\
\cmidrule(lr){4-7} \cmidrule(lr){8-11}
Region & \#Periods & \#Sce & \#Iter & \%Gap & Runtime & Runtime/Iter & \#Iter & \%Gap & Runtime & Runtime/Iter \\
\midrule
SE & 3 & 16 & 17 & 0.02 & 30 & 1.79 & 15 & 0.02 & 29 & 1.91 \\
& 6 & 1,024 & 23 & 0.01 & 171 & 7.44 & 24 & 0.01 & 192 & 8.02 \\
& 9 & 65,536 & 29 & 0.13 & 408 & 14.07 & 18 & 0.13 & 286 & 15.88 \\
& 12 & 4,194,304 & 31 & 0.51 & 717 & 23.12 & 33 & 0.21 & 1,235 & 37.42 \\
EAST & 3 & 16 & 24 & 0.45 & 60 & 2.50 & 26 & 0.45 & 54 & 2.09 \\
& 6 & 1,024 & 98 & 0.12 & 1,515 & 15.46 & 74 & 0.12 & 995 & 13.45 \\
& 9 & 65,536 & 167 & 0.07 & 8,129 & 48.67 & 244 & 0.07 & 21,362 & 87.55 \\
& 12 & 4,194,304 & 134 & 0.54 & 9,705 & 72.42 & 115 & 0.84 & 12,962 & 112.71 \\
US & 3 & 16 & 54 & 0.05 & 190 & 3.52 & 51 & 0.05 & 177 & 3.47 \\
& 6 & 1,024 & 428 & 3.15 & 38,248 & 89.37 & 204 & 0.11 & 8,808 & 43.18 \\
& 9 & 65,536 & 386 & 5.52 & 90,629 & 234.79 & 399 & 3.54 & 90,685 & 227.28 \\
& 12 & 4,194,304 & 250 & 13.21 & 92,013 & 368.05 & 234 & 18.56 & 91,685 & 391.81 \\
\bottomrule
\end{tabular}
}
\end{center}
\vspace*{-6pt}
\end{table}
}

In the second experiment, we evaluate the efficiency and scalability of our enhanced SDDiP framework by comparing four SDDiP variants. Except that the \%Gap column now marks the final optimality gap given by the SDDiP algorithm, the rest of the performance metrics remain consistent with the first experiment.

In Table \ref{tab:exp1.2.1}, the SDDiP versions with the alternating cut strategy achieve the best overall performance across instances. Comparing (parallel) SDDiP with and without the alternating cut strategy, we find that the inclusion of the alternating cuts consistently accelerates convergence and delivers significant performance improvements. For large-scale instances, such as the US region with 9- and 12-period horizons, the alternating cut strategy reduces the optimality gap by nearly half or more relative to the corresponding version without it. On the other hand, the parallel version shows runtime improvements over the sequential one in small- and medium-scale cases (e.g., the SE region and shorter planning horizons in the EAST and US regions). For example, the parallel SDDiP with the alternating cut strategy reduces runtime by one-third in the EAST region with 6 periods compared to the non-parallel version, while maintaining the same optimality gap. However, for larger-scale problems (e.g., the US region with a 12-period planning horizon), we can see that the advantage of parallelization diminishes. Overall, the parallel SDDiP with the alternating cut strategy delivers the best performance across most two-capacity-level instances.

\subsection{Value-of-Stochastic-Solution Metrics for Multi-Stage Stochastic Optimization}
\label{appendix_exp3}

This section investigates the advantages of explicitly accounting for uncertainty through dynamic stochastic modeling, compared to its deterministic counterpart. We adopt the value of the stochastic solution (VSS) metric, which is originally proposed for two-stage stochastic optimizations \cite{birge1982value} and extended to multi-stage settings by \cite{escudero2007value}. In our context, we compute the related metrics as follows:

{
\begin{itemize}
    \item The expected value (EV) denotes the optimal value obtained by solving the deterministic counterpart of the original model (\ref{eq:51})-(\ref{eq:517}) as outlined in Appendix \ref{appendix_c}, which does not consider any disruptions, and the uncertain demand is replaced by the expected value of the average demand scenario.
    
    \item The expectation of expected value at stage $T$, denoted as EEV$_T$, is computed by first fixing the capacity level decisions $Y$ from stage $1$ to stage $T-1$ to the optimal capacity level decisions obtained from the EV solution $\bar{Y}$.

    \begin{equation}
    EEV_t =
    \begin{cases}
        \text{stochastic model }(\ref{eq:1})-(\ref{eq:17}) \\
        \text{s.t.} \quad Y_1(\omega) = \bar{Y}_1 \quad \forall \boldsymbol{\omega} \in \boldsymbol{\Omega}, \\
        \quad\quad \vdots \\
        \quad\quad Y_{t-1}(\omega) = \bar{Y}_{t-1} \quad \forall \boldsymbol{\omega} \in \boldsymbol{\Omega}.
    \end{cases}
    \end{equation}
    
    \hspace{-0.2in} Given these fixed decisions, the model is solved over the full scenario set with the remaining variables adjusting optimally, and $EEV_t$ is defined as the resulting optimal value.
    
    \item The recourse problem (RP) value represents the optimal value attained by directly solving the stochastic model (\ref{eq:1})-(\ref{eq:17}) on the full scenario set.
    
    \item The value of the stochastic solution at stage $T$, denoted as VSS$_T$, serves as a metric to quantify the benefit of explicitly modeling uncertainty through stochastic programming rather than using deterministic approximations. It is defined as the difference between the EEV$_T$ and RP, i.e., VSS$_T = $EEV$_T - $RP.
\end{itemize}
}

\subsection{Complementary Experimental Results for Sensitivity Analysis}
\label{appendix_exp4}

Table \ref{tab:exp4} summarizes the results of the cost parameter sensitivity analysis conducted under different combinations of capacity, relocation, and outsourcing cost levels. For each scenario, the table reports the optimal total costs obtained under three design configurations: non-modular and non-mobile, modular and non-mobile, and modular and mobile capacity. The last three columns list the corresponding performance metrics, including value of modularity (VMoD), value of mobility (VMoB), and value of modular mobility (VMM), which quantify the economic benefits of modularity and mobility. 

\begin{table}[h!]
\begin{center}
\caption{{Cost parameter sensitivity analysis under different modularity and mobility options.}}
\label{tab:exp4}
\resizebox{0.98\textwidth}{!}{%
\begin{tabular}{ccccccccrrr}
\toprule
\#Periods & Capacity & Relocation & Outsource &
$z_{\text{Non-mod\& }}^*$ & 
$z_{\text{Mod\& }}^*$ & 
$z_{\text{Mod\& }}^*$ & 
VMoD & VMoB & VMM \\
  & cost level & cost level & cost level &
$_{\text{Non-mob}}$ & 
$_{\text{Non-mob}}$ & 
$_{\text{Mob}}$ & 
  &   &  \\
\midrule
3 & M & M & M & 5,090,949 & 4,325,071 & 4,233,056 & 765,878 & 92,015 & 857,893 \\
  & M & M & L & 4,114,792 & 3,537,138 & 3,455,307 & 577,654 & 81,831 & 659,485 \\
  & M & M & H & 6,529,413 & 5,643,085 & 5,523,885 & 886,328 & 119,200 & 1,005,528 \\
  & M & L & M & 5,090,949 & 4,325,071 & 4,228,857 & 765,878 & 96,214 & 862,092 \\
  & M & H & M & 5,090,949 & 4,325,071 & 4,241,311 & 765,878 & 83,760 & 849,638 \\
  & L & M & M & 3,893,174 & 3,384,763 & 3,271,947 & 508,411 & 112,816 & 621,227 \\
  & H & M & M & 7,217,814 & 6,043,177 & 5,945,855 & 1,174,637 & 97,322 & 1,271,959 \\
\vspace*{-6pt}\\
6 & M & M & M & 6,019,385 & 5,332,299 & 5,208,949 & 687,086 & 123,350 & 810,436 \\
  & M & M & L & 5,482,602 & 4,814,916 & 4,696,075 & 667,686 & 118,841 & 786,527 \\
  & M & M & H & 7,002,063 & 5,958,858 & 5,801,574 & 1,043,205 & 157,284 & 1,200,489 \\
  & M & L & M & 6,019,385 & 5,332,299 & 5,206,307 & 687,086 & 125,992 & 813,078 \\
  & M & H & M & 6,019,385 & 5,332,299 & 5,214,234 & 687,086 & 118,065 & 805,151 \\
  & L & M & M & 4,262,542 & 3,742,180 & 3,624,817 & 520,362 & 117,363 & 637,725 \\
  & H & M & M & 9,537,997 & 8,274,569 & 8,118,520 & 1,263,428 & 156,049 & 1,419,477 \\
\bottomrule
\end{tabular}
}
\vspace{-0.2in}
\end{center}
\end{table}


\section{Deterministic Model for the DSMMCP Problem}
\label{appendix_c}

In the deterministic model, the throughput rate per capacity module is fixed at its maximum level, implying no disruption. Moreover, the average demand across scenarios is utilized as the given demand parameters. The formulation is presented as follows:

{
\renewcommand{\baselinestretch}{0.75}\selectfont
\begin{align}
    &\min \quad
    \sum_{t \in \boldsymbol{T}} \Bigg(
        \sum_{j \in \boldsymbol{J}} 
        \sum_{l_1 \in \boldsymbol{L}} 
        \sum_{l_2 \in \boldsymbol{L}} 
        c_{j l_1 l_2 t}^F  Y_{j l_1 l_2 t} \nonumber \\
    & \qquad
        + \sum_{(j,j^{'}) \in \boldsymbol{\mathit{JJ}}} 
        c_{j j^{'} t}^A  F_{jj^{'}t}  \nonumber\\
    & \qquad
        + \sum_{(i,j) \in \boldsymbol{\mathit{IJ}}} 
        c_{i j t}^T  X_{ijt} 
        + \sum_{j \in \boldsymbol{J}} 
        c_{t}^O  R_{jt}
    \Bigg)
    & \label{eq:51}\\
    & \text{s.t.,} \nonumber\\
    & \sum_{j \in \boldsymbol{J} \mid (i,j) \in \boldsymbol{\mathit{IJ}}} 
    X_{ijt} \geq d_{it}, \quad \forall i \in \boldsymbol{I}, t \in \boldsymbol{T}
    & \label{eq:52} \\
    & k_{jt} S_{jt} + R_{jt}
    \geq 
    \sum_{i \in \boldsymbol{I} \mid (i,j) \in \boldsymbol{\mathit{IJ}}} 
    X_{ijt}, \nonumber\\
    & \qquad \qquad \qquad \qquad \qquad \qquad \forall j \in \boldsymbol{J}, t \in \boldsymbol{T}
    & \label{eq:54} \\
    & S_{j1}
    = u_{j v_j}
    + \sum_{j^{'} \in \boldsymbol{J} \cup \{0\} \mid (j,j^{'}) \in \boldsymbol{\mathit{JJ}}}
    \left(
        F_{j^{'}j1} - F_{jj^{'}1}
    \right), \nonumber\\
    & \qquad \qquad \qquad \qquad \qquad \qquad \qquad \forall j \in \boldsymbol{J}
    & \label{eq:54.5} \\
    & S_{jt}
    = S_{j,t-1}
    + \sum_{j^{'} \in \boldsymbol{J} \cup \{0\} \mid (j,j^{'}) \in \boldsymbol{\mathit{JJ}}}
    \left(
        F_{j^{'}jt} - F_{jj^{'}t}
    \right), \nonumber\\
    & \qquad \qquad \qquad \qquad \qquad \quad \forall j \in \boldsymbol{J},\ t \in \boldsymbol{T} \backslash \{1\}
    & \label{eq:55} \\
    & S_{jt}
    = \sum_{l_1 \in \boldsymbol{L}}
      \sum_{l_2 \in \boldsymbol{L}}
      u_{j l_2} \, Y_{j l_1 l_2 t}, \quad \forall j \in \boldsymbol{J}, t \in \boldsymbol{T}
    & \label{eq:56} \\
    & \sum_{l \in \boldsymbol{L}}
    Y_{j v_j l 1} = 1, \quad \forall j \in \boldsymbol{J}
    & \label{eq:57} \\
    & \sum_{l_1 \in \boldsymbol{L}}
    Y_{j, l_1, l, t-1}
    =
    \sum_{l_2 \in \boldsymbol{L}}
    Y_{j, l, l_2, t}, \nonumber\\
    & \qquad \qquad \qquad \qquad \forall j \in \boldsymbol{J},\ l \in \boldsymbol{L},\ t \in \boldsymbol{T} \backslash \{1\}
    & \label{eq:58} \\
    & \sum_{l_1 \in \boldsymbol{L}}
      \sum_{l_2 \in \boldsymbol{L}}
      Y_{j l_1 l_2 t}
    = 1, \quad \forall j \in \boldsymbol{J}, t \in \boldsymbol{T}
    & \label{eq:59}\\
    & Y_{j l_1 l_2 t} \leq w_t, \quad \forall j \in \boldsymbol{J},\ l_1 \neq l_2,\ t \in \boldsymbol{T_t}
    & \label{eq:510} \\
    & Y_{j l_1 l_2 t} \in \{0,1\}, \quad \forall j \in \boldsymbol{J},\ l_1,l_2 \in \boldsymbol{L},\ t \in \boldsymbol{T}
    & \label{eq:515} \\
    & F_{jj^{'}t},\ S_{jt} \in \mathbb{Z}_{\ge 0}, \quad \forall (j,j^{'}) \in \boldsymbol{\mathit{JJ}},\ j \in \boldsymbol{J},\ t \in \boldsymbol{T}
    & \label{eq:516} \\
    & X_{ijt},\ R_{jt} \geq 0, \quad \forall (i,j) \in \boldsymbol{\mathit{IJ}},\ j \in \boldsymbol{J},\ t \in \boldsymbol{T}
    & \label{eq:517}.
\end{align}
}

\section{Theoretical Properties of the SDDiP Enhancements}
\label{appendix:theory}

\subsection{Monotone Value of Partial Adaptivity}
\label{appendix:pamssp_proof}

\textbf{Proposition 3.1 }(Monotone Value of Partial Adaptivity).
Let $z^{\mathsf{MSSP}}$ denote the optimal value of the fully adaptive multistage model,
$z^{\mathsf{TSSP}}$ the optimal value of the two-stage static model,
and $z^{\mathsf{PAMSSP}}(a)$ the optimal value of the partially adaptive model with at
most $a$ revision points. Then for any integers $1 \le a_1 \le a_2 \le T$,

{
\[
z^{\mathsf{MSSP}}
\;\le\;
z^{\mathsf{PAMSSP}}(a_2)
\;\le\;
z^{\mathsf{PAMSSP}}(a_1)
\;\le\;
z^{\mathsf{TSSP}},
\]
and equivalently,
\[
z^{\mathsf{PAMSSP}}(a) - z^{\mathsf{PAMSSP}}(a+1) \ge 0 \quad \forall\, a < T.
\]
}

\begin{proof}{Proof.}
Let $\boldsymbol{T}=\{1,\ldots,T\}$ denote the planning horizon, and let
$\mathcal{T}_t$ be the set of nodes at stage $t$ in the scenario tree.
A revision schedule is a vector $\boldsymbol{w}=(w_t)_{t\in\boldsymbol{T}} \in \{0,1\}^T$
chosen at the first stage.  When $w_t=0$, constraints \eqref{eq:10} force
$Y_{j l_1 l_2 n} = 0$ for all $l_1 \neq l_2$, and together with
\eqref{eq:8}--\eqref{eq:9}, this implies that all nodes at stage $t$ must inherit the
capacity level of their parent node $a(n)$; i.e., nonanticipativity is enforced within that
stage block. When $w_t=1$, transitions between capacity levels are allowed. Hence, $w_t$
encodes the admissible branching structure of the state variable.

For a fixed revision schedule $\boldsymbol{w}$, let
\begin{align*}
\mathcal{F}(\boldsymbol{w})
&=
\Big\{
(Y,F,X,R,S) \;\big|\;
\eqref{eq:2}-\eqref{eq:9},\ \eqref{eq:10},\\
& \qquad \eqref{eq:15}-\eqref{eq:17}
\text{ hold under schedule } \boldsymbol{w}
\Big\}
\end{align*}

be the feasible region.  The set of all schedules that allow at most $a$ revisions is
\[
\mathcal{W}(a)
=
\left\{
\boldsymbol{w}\in\{0,1\}^T:\ \sum_{t\in\boldsymbol{T}} w_t \le a,\; w_1 = 1
\right\},
\]
where $w_1=1$ is imposed because the first stage must allow adjustment from the initial
capacity level.

The feasible region of PAMSSP with at most $a$ revisions is the union
\begin{align*}
\mathcal{F}^{\mathsf{PAMSSP}}(a)
=
\bigcup_{\boldsymbol{w}\in\mathcal{W}(a)} \mathcal{F}(\boldsymbol{w}).
\end{align*}

The special cases satisfy
\begin{align*}
\mathcal{F}^{\mathsf{TSSP}} = \mathcal{F}^{\mathsf{PAMSSP}}(1), \
\mathcal{F}^{\mathsf{MSSP}} = \mathcal{F}^{\mathsf{PAMSSP}}(T).
\end{align*}

Since $\mathcal{W}(a_1)\subseteq \mathcal{W}(a_2)$ whenever $a_1\le a_2$, we obtain
the nesting
\begin{align*}
\mathcal{F}^{\mathsf{TSSP}}
\subseteq
\mathcal{F}^{\mathsf{PAMSSP}}(a_1)
\subseteq
\mathcal{F}^{\mathsf{PAMSSP}}(a_2)
\subseteq
\mathcal{F}^{\mathsf{MSSP}}. \hspace{-0.1in}
\end{align*}

All models minimize the same objective function over these nested feasible regions.
Therefore,
\begin{align*}
z^{\mathsf{MSSP}}
\;\le\;
z^{\mathsf{PAMSSP}}(a_2)
\;\le\;
z^{\mathsf{PAMSSP}}(a_1)
\;\le\;
z^{\mathsf{TSSP}}.
\end{align*}

Setting $a_1=a$ and $a_2=a+1$ yields the marginal statement
$z^{\mathsf{PAMSSP}}(a)-z^{\mathsf{PAMSSP}}(a+1)\ge 0$ for all $a<T$.
\end{proof}

\qedhere
\begin{remark}[Value-of-information interpretation]
This proposition is consistent with the classic value-of-information principle in multistage stochastic programming
\citep{birge1997introduction,shapiro2021lectures}: relaxing nonanticipativity constraints -- by introducing additional revision points -- enlarges the feasible region and cannot
worsen the optimal value.  Our PAMSSP formulation operationalizes this principle through
the revision schedule $\boldsymbol{w}$, connecting the general value-of-information logic to the binary-state
SDDiP framework of \cite{zou2018multistage,zou2019stochastic}.
\end{remark}

\subsection{Dominance of Strengthened Pareto-Optimal Cuts}
\label{appendix:strengthened-benders}

\textbf{Proposition 4.1 }(Dominance of Strengthened Pareto-Optimal Cuts).
Let $\hat{Y}_n^k$ and $Y^{\circ k}_n$ denote the state and core point selected at iteration $k$ during the forward pass. For each child node $\nu \in \mathcal{C}(n)$, let $\rho_\nu^k$ and $\alpha_\nu^k$ denote the optimal value and optimal dual multipliers of the Magnanti-Wong secondary problem $M_\nu^k(\hat{Y}_n^k, Y^{\circ k}_n)$ enforcing $Y^P_\nu = \hat{Y}_n^k$.
Let $\eta_\nu^k$ be the optimal value of the associated Lagrangian relaxation constructed using the same multipliers $\alpha_\nu^k$.
Then, for each feasible parent-state vector $Y_n$,
\begin{align*}
 &\sum_{\nu \in \boldsymbol{\mathit{C}(n)}} p_{n \nu} \Big( \eta_{\nu}^{k} + (\alpha^{k}_{\nu})^{\intercal}  Y_n \Big)
\\
& \quad \;\ge\;
\sum_{\nu \in \boldsymbol{\mathit{C}(n)}} p_{n \nu} \Big( \rho_{\nu}^{k} + (\alpha^{k}_{\nu})^{\intercal} (Y_n - Y^{\circ k}_n) \Big).
\end{align*}

Hence, the strengthened Pareto-optimal cut is never weaker than the classical Pareto-optimal cut.

\vspace{0.2in}
\begin{proof}{Proof.}
For a fixed child node $\nu$, LP relaxation of the forward subproblem $R^k_\nu$ can be written as
\begin{flalign}
\sigma_\nu^k
=
\min_{(x,Y^P_\nu)\in \mathcal{F}_{LP}}
\left\{\, f_\nu(x)
    + (\pi_\nu^k)^\top(\hat{Y}_n^k - Y^P_\nu)
\right\}, \label{eq:66}
\end{flalign}

where $\mathcal{F}_{LP}$ is the LP-relaxed feasible region. 
By strong duality, $\pi_\nu^k$ is the maximizing dual vector, so substituting it yields an exact
representation of $\sigma_\nu^k$.

Moreover, we assume the dual of $R^k_\nu$ can be written as:
\begin{flalign*}
D^k_\nu(\hat{Y}^k_n) 
:=\;& \max
    g_\nu(\pi_\nu) + (\hat{Y}^k_n)^{\intercal} \lambda_\nu - \bm{1}^{\intercal} \nu_\nu \\
& \text{s.t. } (\pi_\nu, \lambda_\nu, \mu_\nu, \nu_\nu)
    \in \boldsymbol{\Lambda_\nu},
\end{flalign*}

where $g_\nu(\cdot)$ represents the dual objective function, while $\pi_\nu$ is dual variable associated with these internal constraints. In addition, dual variable $\lambda_\nu$ corresponds to the linking equality, and $\mu_\nu, \nu_\nu \ge 0$ correspond to the bound constraints $0 \le Y^P_\nu \le 1$. $\boldsymbol{\Lambda_\nu}$ denotes the feasible set of dual multipliers. 

Then, we construct the Magnanti-Wong secondary problem:
\begin{flalign}
    M^k_\nu(\hat{Y}^k_n, Y^{\circ k}_n) :=\;& \max
    g_\nu(\pi_\nu) + (Y^{\circ k}_n)^{\intercal} \lambda_\nu - \bm{1}^{\intercal} \nu_\nu \label{eq:dual_compact}\\
    & \text{s.t. } (\pi_\nu, \lambda_\nu, \mu_\nu, \nu_\nu)
    \in \boldsymbol{\Lambda_\nu}\\
    & g_\nu(\pi_\nu) + (\hat{Y}^k_n)^{\intercal} \lambda_\nu - \bm{1}^{\intercal} \nu_\nu = \sigma_\nu^{k},\label{eq:65}
\end{flalign}

where $Y^{\circ k}_n$ is the core point. Solving $M^k_\nu$ yields the optimal value $\rho_\nu^{k}$ and the optimal solution of $\lambda_\nu$, denoted by $\alpha_\nu^{k}$ (i.e., $\alpha_\nu^{k}:= \lambda^*_\nu$). Moreover, it is important to note that $\alpha_\nu^{k}$ also serves as the optimal dual multiplier associated with the linking constraint of $R^k_\nu$. Since $M^k_\nu$ enforces constraint \eqref{eq:65} under dual feasibility, weak duality implies that this equality can hold only at dual-optimal solutions of $R^k_\nu$. Therefore, $\alpha_\nu^{k}$ lies on the dual-optimal face of the dual of $R^k_\nu$.

Therefore, using the Magnanti-Wong optimality constraint (\ref{eq:66}) and (\ref{eq:65}), we have:
\begin{align*}
&\rho_\nu^{k}= g_\nu(\pi_\nu^*)+(Y_n^{\circ k})^\top \alpha_\nu^{k}-\bm 1^\top \nu_\nu^* \\
&= \bigl[g_\nu(\pi_\nu^*)+(\hat Y_n^{k})^\top \alpha_\nu^{k}-\bm 1^\top \nu_\nu^*\bigr]
+(\alpha_\nu^{k})^\top(Y_n^{\circ k}-\hat Y_n^{k}) \\
&= \sigma_\nu^{k}+(\alpha_\nu^{k})^\top(Y_n^{\circ k}-\hat Y_n^{k})\\
&=\min_{(x,Y^P_\nu)\in \mathcal{F}_{LP}}
\left\{
    f_\nu(x) + (\alpha_\nu^k)^\top(Y_n^{\circ k} - Y^P_\nu)
\right\}.
\end{align*}

On the other hand, the strengthened cut evaluates the same Lagrangian function over the mixed-integer
feasible region $\mathcal{F}_{INT} \subseteq \mathcal{F}_{LP}$:
\begin{align*}
&\eta_\nu^k + (\alpha_\nu^k)^\top Y_n^{\circ k}\\
& \quad = \min_{(x,Y^P_\nu)\in \mathcal{F}_{INT}}
\left\{
    f_\nu(x) + (\alpha_\nu^k)^\top(Y_n^{\circ k} - Y^P_\nu)
\right\}.
\end{align*}

Since the minimization is taken over a subset of $\mathcal{F}_{LP}$,
\[
\eta_\nu^k + (\alpha_\nu^k)^\top Y_n^{\circ k}
\;\ge\;
\rho_\nu^k.
\]

Rearranging and adding $(\alpha_\nu^k)^\top Y_n$ to both sides yields
\[
\eta_\nu^k + (\alpha_\nu^k)^\top Y_n
\;\ge\;
\rho_\nu^k + (\alpha_\nu^k)^\top (Y_n - Y_n^{\circ k}).
\]

Finally, summing with probabilities $p_{n\nu}$ over all child nodes preserves the inequality,
which proves the dominance of the strengthened cut.
\end{proof}

\begin{remark}[Strict dominance]
Strict dominance holds whenever the LP relaxation produces fractional values of $Y_\nu$
for at least one child node.  
In this case, the mixed-integer feasible region used in the strengthened cut excludes the
LP-optimal solution, leading to
$\eta_\nu^k + ({\alpha_\nu^{k}})^{\top}Y^{\circ k}_n > \rho_\nu^k$.
Hence, the strengthened cut strictly dominates the classical cut at $Y_n = Y_n^{\circ k}$.
\end{remark}

\begin{remark}[Geometric insight]
Both cuts share the same slope vector $\sum_{\nu} p_{n\nu}\alpha_\nu^k$, but the strengthened
cut has a higher intercept because $\eta_\nu^k \ge
\rho_\nu^{k} - ({\alpha_\nu^{k}})^{\top}Y^{\circ k}_n$ for each child node $\nu$.
Thus, the strengthened cut is a vertical upward shift of the classical cut and dominates it
pointwise.  
This is consistent with the Pareto-optimal cut theory of
\cite{magnanti1981accelerating,papadakos2008practical} and with strengthened SDDiP cuts
in \cite{zou2019stochastic}.
\end{remark}

\end{appendices}

\end{document}